\documentclass[conference]{IEEEtran}
\usepackage[utf8]{inputenc}

\usepackage{amsmath}
\usepackage{amsfonts}
\usepackage{amssymb}

\usepackage{times}
\usepackage{mathtools}
\usepackage{calrsfs}
\usepackage{comment}
\usepackage{float}
\usepackage{graphicx}
\usepackage{caption}
\usepackage{subcaption}
\usepackage{makecell, cellspace} 
\usepackage{siunitx,booktabs}
\DeclareMathAlphabet{\pazocal}{OMS}{zplm}{m}{n}

\DeclareMathOperator*{\argmin}{arg\,min}\usepackage[ruled,longend]{algorithm2e}

\usepackage{algorithmic}
\usepackage{array}
\usepackage{bm}
\usepackage{authblk}
\newcommand*{\affaddr}[1]{#1} 
\newcommand*{\affmark}[1][*]{\textsuperscript{#1}}
\newcommand*{\email}[1]{\texttt{#1}}

\usepackage[numbers]{natbib}

\usepackage{multicol}
\usepackage[bookmarks=true]{hyperref}

\newcommand{\tr}[1]{{#1}^\mathsf{T}}

\usepackage{color}
\usepackage{xcolor}

\begin{document}

\title{Constrained Differential Dynamic Programming  Revisited}



\author{
Yuichiro Aoyama\affmark[1,2], George Boutselis\affmark[1], Akash Patel\affmark[1], and Evangelos A. Theodorou\affmark[1]\\
\affaddr{\affmark[1]School of Aerospace Engineering, Georgia Institute of Technology, Atlanta, GA, USA}\\
\affaddr{\affmark[2]Komatsu Ltd., Tokyo, Japan}\\
\email{\{yaoyama3, apatel435, gbouts, evangelos.theodorou\}@gatech.edu}\\
}



%

\maketitle

\begin{abstract}
Differential Dynamic Programming (DDP) has become a well established method for unconstrained trajectory optimization. Despite its several applications in robotics and controls however, a widely successful constrained version of the algorithm has yet to be developed. This paper builds upon penalty methods and active-set approaches, towards designing a Dynamic Programming-based methodology for constrained optimal control. Regarding the former, our derivation employs a constrained version of Bellman's principle of optimality, by introducing a set of auxiliary slack variables in the backward pass. In parallel, we show how Augmented Lagrangian methods can be naturally incorporated within DDP, by utilizing a particular set of penalty-Lagrangian functions that preserve second-order differentiability. We demonstrate experimentally that our extensions (individually and combinations thereof) enhance significantly the convergence properties of the algorithm, and outperform previous approaches on a large number of simulated scenarios.
\end{abstract}

\IEEEpeerreviewmaketitle

\section{Introduction}

Trajectory optimization problems arise very frequently in robotics and controls applications. Examples include finding suitable motions for robotic grasping and manipulation tasks \cite{kumar2016optimal}, or minimizing fuel for orbital transfers \cite{lantoine2008hybrid}. Mathematically speaking, such problems require computing a state/control sequence that minimizes a specified cost function, while satisfying the dynamics constraints of the agent. Common methodologies for trajectory optimization rely on optimal control and/or optimization theory. The former approach provides fundamental principles for obtaining solutions (based, for example, on Dynamic Programming or the Hamilton-Jacobi-Bellman equation), which, however, do not scale well with high-dimensional, nonlinear problems \cite{bertsekas1995dynamic}. In contrast, standard direct optimization methods can be used for discrete optimal control \cite{mombaur2009using}. The main drawback of these works is that feasibility with respect to dynamics has to be explicitly imposed, thus slowing down the optimization process \cite{xie2017differential}.

One of the most successful trajectory optimization algorithms is Differential Dynamic Programming (DDP), originally developed by Jacobson and Mayne \cite{ddp}. DDP is an indirect method which utilizes Bellman's principle of optimality to split the problem into ``smaller" optimization subproblems at each time step. Under mild assumptions on the cost and dynamics, it can be shown that DDP achieves locally quadratic convergence rates \cite{liao1991convergence}. While the original method relies on second-order derivatives, one of its variations, iterative-Linear-Quadratic-Regulator (iLQR), uses only Gauss-Newton approximations of the cost Hessians as well as first-order expansions of the dynamics \cite{li2004iterative}, which is often numerically advantageous. The aforementioned algorithms have been employed in various applications such as robotic manipulation \cite{kumar2016optimal}, bipedal walking \cite{morimoto2003minimax} and model-based reinforcement learning \cite{levine2014learning}, to name a few.

While unconstrained DDP has been widely tested and used over the past decades, its constrained counterpart has yet to be properly established. Since most practical applications in controls and robotics include state and/or control constraints (e.g., navigating through obstacles, respecting joint/actuator limits, etc.), off-the-shelf optimization solvers still remain the most popular tool for trajectory optimization among scientists and practitioners \cite{kobilarov2011discrete, mombaur2009using}. A few works have attempted to extend the DDP framework to the constrained case. \cite{tassa2014control, murray1979constrained} considered the case of control bounds by solving several quadratic programs over the trajectory, and \cite{giftthaler2017projection} dealt with equality constraints only via projection techniques. The works in \cite{lantoine2008hybrid, xie2017differential} utilized the Karush-Kuhn-Tucker (KKT) conditions when both state and control constraints are present, with \cite{xie2017differential}, \cite{Lin91}, in particular, solving successive quadratic programs in the forward pass of the method. \cite{lantoine2008hybrid, plancher2017constrained}, \cite{howell2019altro} also discussed combining DDP with an Augmented Lagrangian (AL) approach; \cite{lantoine2008hybrid} updated the Lagrange multipliers via forward/backward passes, while \cite{plancher2017constrained}, \cite{howell2019altro} utilized schemes from the standard Powell-Hestenes-Rockafellar (PHR) methodology \cite{nocedal2006numerical} with first-order approximations of the Hessians.

In this paper we build upon the works in \cite{lantoine2008hybrid, xie2017differential, plancher2017constrained} to develop a state- and control-constrained version of DDP in discrete time. Specifically, we extend \cite{lantoine2008hybrid, xie2017differential} by introducing a slack variable formulation into Bellman's principle, and thus avoid assumptions regarding the active constraints of the problem. Moreover, we propose an Augmented Lagrangian-inspired algorithm, by considering a set of penalty functions that preserves smoothness of the transformed objective function. This property was not satisfied in \cite{plancher2017constrained}, but is required to establish the convergence properties of DDP \cite{liao1991convergence}. These two methodologies can be used separately, or be properly combined for improved numerical performance.

We will save the in-depth discussion about technical differences between our methods and previous papers for subsequent sections. Nevertheless, we note that a comparison among different constrained optimization methods on various simulated scenarios will be provided, which will highlight the efficiency and generalizability of our approach; something which has been lacking from previous DDP-related schemes. To the best of the authors' knowledge, such an extensive experimental study on constrained trajectory optimization has not been conducted in the past. We believe that the current work is a key step towards the development of a numerically robust, constrained version of Differential Dynamic Programming, and opens up multiple directions for research and further improvements.

The remaining of this paper is organized as follows: Section \ref{sec:preliminaries} gives some basics on unconstrained DDP and constrained optimization theory. In Section \ref{sec:CDDP} we explain our KKT-based DDP algorithm with slack variables (S-KKT), while Section \ref{sec:AL-DDP} discusses our AL-inspired method, as well as a combination thereof. Numerical experiments and in-depth comparisons between our methodologies and previous implementations are provided in Section \ref{sec:results}. Section \ref{sec:conclusion} concludes the paper and discusses possible extensions of the current work.


\section{Preliminaries}
\label{sec:preliminaries}

\subsection{Unconstrained Differential Dynamic Programming}
\label{subsec:ddp}
We will briefly cover here the derivation and implementation of Differential Dynamic Programming (DDP). More details can be found in \cite{ddp, li2004iterative}.

Consider the discrete-time optimal control problem
\begin{equation}
\label{unconstrained_optimalcontrol}
\begin{split}
\min_{\bm{U}}&\hspace{0.8mm}J(\bm{X},\bm{U})=\min_{\bm{U}}\big[{\textstyle\sum_{k=0}^{N-1} l(\bm{x}_k,\bm{u}_k)}+\phi(\bm{x}_N)\big]\\
&\text{subject to:}\hspace{3mm} \bm{x}_{k+1}=\bm{f}(\bm{x}_k,\bm{u}_k),\quad k=0,...,N-1.
\end{split}
\end{equation}
where $\bm{x}_k\in\mathbb{R}^n$, $\bm{u}_k\in\mathbb{R}^m$ denote the state and control input of the system at time instant $t_k$, respectively, and $\bm{f}:\mathbb{R}^n\times \mathbb{R}^m\rightarrow\mathbb{R}^n$ corresponds to the transition dynamics function. The scalar-valued functions $l(\cdot,\cdot)$, $\phi(\cdot)$,  $J(\cdot)$ denote the running, terminal and total cost of the problem, respectively. We also let $\bm{X} :=(\tr{\bm{x}}_{0},\dots,\tr{\bm{x}}_{N})$, $\bm{U} :=(\tr{\bm{u}}_{0},\dots,\tr{\bm{u}}_{N-1})$ be the state/control sequences over the horizon $N$.

Of paramount importance is the concept of the {\it value function}, which represents the minimum cost-to-go at each state and time. It is defined as:
\begin{equation}
\label{value-function-definition}
V_k(\bm{x}_k):=\min_{\bm{u}_k}J(\bm{X},\bm{U}).
\end{equation}
Based on this, {\it Bellman's principle of optimality} gives the following rule:
\begin{equation}\label{bellman}
    V_k(\bm{x}_k) =  \min_{\bm{u}_k} [l(\bm{x}_k,\bm{u}_k) + V_{k+1}(\bm{x}_{k+1})].
\end{equation}
DDP finds local solutions to \eqref{unconstrained_optimalcontrol} by expanding both sides of \eqref{bellman} about given nominal trajectories, $\bar{\bm{ X}}$, $\bar{\bm{ U}}$. Specifically, let us define the $Q$ function as the argument of $\min$ on the right-hand-side of \eqref{bellman}:
\begin{equation}
\label{Qfunction}
Q_k(\bm{x}_k,\bm{u}_k) = l(\bm{x}_k,\bm{u}_k) + V_{k+1}(\bm{x}_{k+1}).
\end{equation}
We now proceed by taking quadratic expansions of $Q_k$ about $\bar{\bm{ X}}$, $\bar{\bm{ U}}$. According to \eqref{Qfunction} and an additional expansion of $\bm{x}_{k+1}=\bm{f}(\bm{x}_{k},\bm{u}_{k})$, this will give
\begin{equation*}
    \begin{split}
  &Q_k(\bm{x}_k,\bm{u}_k)\approx  Q_k+{Q}_{\bm{x},k}^{\mathsf{T}}\delta\bm{x}_k+{Q}_{\bm{u},k}^{\mathsf{T}}\delta\bm{u}_k+\\
  &\quad\quad\textstyle{\frac{1}{2}}(\tr{\delta\bm{x}}_k{Q}_{\bm{xx},k}\delta\bm{x}_k+2\tr{\delta\bm{x}}_k{Q}_{\bm{xu},k}\delta\bm{u}_k+\tr{\delta\bm{u}}_k{Q}_{\bm{uu},k}\delta\bm{u}_k),
  \end{split}
\end{equation*}
with
\begin{equation}
\label{Qexpanded}
\begin{split}
    &{Q}_{\bm{xx},k} ={l}_{\bm{xx}}+\bm{f_x}^{\mathsf{T}}{V}_{\bm{xx},k+1}\bm{f_x},\hspace{0.8mm}{Q}_{\bm{x},k}={l}_{\bm{x}}+\bm{f_x}^{\mathsf{T}}{V}_{\bm{x},k+1}
     \\
    &{Q}_{\bm{uu},k} = {l_{\bm{uu}}}+\bm{f_u}^{\mathsf{T}}{{V}_{\bm{xx},k+1}}\bm{f_u}, \hspace{0.8mm}{Q}_{\bm{u},k}={l_{\bm u}}+\bm{f_u}^{\mathsf{T}}{{V}_{\bm{x},k+1}}
    \\
    &{Q}_{\bm{xu},k} ={l_{\bm{xu}}}+\bm{f_x}^{\mathsf{T}}{{V}_{\bm{xx},k+1}}\bm{f_u}.
\end{split}
\end{equation}
Here $\delta \bm{x}_k:=\bm{x}_k -\bar{\bm{x}}_k$, $\delta\bm{u}_k:=\bm{u}_k - \bar{\bm{u}}_k$ are deviations about the nominal sequences. It is also implied that the $Q$ functions above are evaluated on $\bar{\bm{ X}}$, $\bar{\bm{ U}}$.

After plugging \eqref{Qexpanded} into \eqref{bellman}, we can explicitly optimize with respect to $\delta\bm{u}$ and compute the locally optimal control deviations. These will be given by
\begin{equation}
\label{delta-u-star}
\begin{split}
    & \delta \bm{u}^{\ast}_k=\bm{k}_k+\bm{K}_k\delta \bm{x}_{k},\\
    \text{with}\quad &\bm{k} := -{Q}^{-1}_{\bm{uu}}{Q_{\bm{u}}},\hspace{1.8mm} \bm{K} = -{Q}^{-1}_{\bm{uu}}{Q_{\bm{ux}}}.
\end{split}
\end{equation}
Finally, observe that $\delta \bm{u}^\ast$ requires knowledge of the value function on the nominal rollout. To this end, $V_k$ will be quadratically expanded and will be plugged along with $\delta \bm{u}^\ast_k$ into \eqref{bellman} to give:
\begin{equation} \label{eq:value_riccati}
\begin{split}
    V_{\bm{x},k}&={Q_{\bm{x},k}}-{Q_{\bm{xu},k}}{Q_{\bm{uu},k}^{-1}}{Q_{\bm{u},k}}\\
    V_{\bm{xx},k}&={Q_{\bm{xx},k}}-{Q_{\bm{xu},k}}Q_{\bm{uu},k}^{-1}{Q_{\bm{ux},k}}.
    \end{split}
\end{equation}
The equations above are propagated backwards in time, since at the final horizon the value function equals the terminal cost. After the backward pass is complete, a new state-control sequence is determined in a forward pass, and this trajectory
is then treated as the new nominal trajectory for the next iteration. The procedure is then repeated until certain convergence criteria are satisfied.

To ensure convergence, $Q_{\bm{uu}}$ must be regularized, when its positive definiteness cannot be guaranteed \cite{liao1991convergence}. Typically, line-search on $\delta\bm{u}^\ast$ is also performed with respect to the total  cost  in the forward pass. We finally note that in this paper we consider only first-order expansions of the dynamics as in \cite{li2004iterative}, which tends to be less computationally expensive and more numerically stable than using second-order terms.

\subsection{Constrained optimization theory}
We present here preliminaries on constrained optimization. Due to space limitations, we only consider inequality constraints, though similar results hold for equality constraints.

\subsubsection{KKT conditions}
\label{subsec:kkt}
Consider the optimization problem
\begin{align} \label{constrainedOptimization}
  \min_{\bm{x}} h(\bm{x})&\\\notag
     \text{subject to}\quad \bm{g}(\bm{x})  &\leq \bm{0},
     \end{align}
where $\bm{g} = (g_{1}(\bm{x}),\dots,g_{w}(\bm{x}))^{\mathsf{T}}$ is a vector of $w$ constraints. The {\it Lagrangian} is defined as follows:
\begin{align}
    L = h(\bm{x}) + \bm{\lambda}^{\mathsf{T}}\bm{g}(\bm{x}),
\end{align}
for a real vector $\bm{\lambda}\in\mathbb{R}^w$. The
 {\it Karush–Kuhn–Tucker} (KKT) conditions are necessary optimality conditions for problems of the type \eqref{constrainedOptimization}, and state that a local solution must satisfy \cite[Section 12]{nocedal2006numerical}
\begin{align}\label{eq:first_optimality}
\notag\frac{\partial L}{\partial \bm{x}} = \frac{\partial h(\bm{x})}{\partial \bm{x}} +\frac{\partial \bm{g}(\bm{x})}{\partial \bm{x}}^{\mathsf{T}}\bm{\lambda} = 0\\
g_{i}(\bm{x})\leq0,\quad{\lambda_i} \geq {0}, \hspace{3mm}\text{and}\hspace{3mm} {\lambda} _{i} g_{i}(\bm{x}) ={0},\hspace{1.7mm}\forall i,
\end{align}
where $\lambda$ now corresponds to the {\it Lagrange multipliers}. For this result, we need to assume differentiability of $f$ and $\bm{g}$, as well as linear independence of the gradients of active constraints.

\subsubsection{Augmented Lagrangian} \label{subsec:AL}
Let us define the Augmented Lagrangian function corresponding to \eqref{constrainedOptimization} as \cite{birgin2005numerical, nocedal2006numerical}

\begin{equation}
\label{augmentedLagrangian}
L_A(\bm{x}, \bm{\lambda}, \bm{\mu}):=h(\bm{x})+\textstyle\sum_{i}\mathcal{P}(g_i(\bm{x}), \lambda_{i}, \mu_{i}),
\end{equation}
where $\bm\lambda$, $\bm\mu$ correspond to the Lagrange multipliers and penalty parameters respectively, while $\mathcal{P}(\cdot)$ is the penalty function for inequalities. When $\mathcal{P}(\cdot)$ satisfies certain properties, it can be shown that minimization of \eqref{augmentedLagrangian} can give a solution to \eqref{constrainedOptimization}, under mild assumptions \cite{birgin2005numerical}.

Loosely speaking, the corresponding optimization process can be divided into an inner and outer loop. In the inner loop, a local minimizer is found for \eqref{augmentedLagrangian} by an unconstrained optimization methodology. At the outer loop, the Lagrange multipliers are updated as: $\lambda_{i }\leftarrow \mathcal{P}'(g_{i},\lambda_{i},\mu_{i})$, where $\mathcal{P}'(y, \lambda, \mu):=\frac{\partial}{\partial y}\mathcal{P}(y,\lambda,\mu)$. Moreover, the penalty parameters are increased monotonically, when constraint improvement is not satisfactory.

The most popular Augmented Lagrangian algorithm uses the penalty function $P(y,\lambda,\mu)=\frac{1}{2\mu}(\max(0,\lambda+\mu y)^2-\lambda^2)$ and is known as the Powell-Hestenes-Rockafellar (PHR) method \cite{nocedal2006numerical}. Despite its success, one key drawback is that the objective function of each subproblem is not twice differentiable, which may cause numerical instabilities when used within second-order algorithms \cite{birgin2005numerical}.

For completeness, we give the required properties for $\mathcal{P}$ in the appendix, and refer the interested reader to \cite[Section 17]{nocedal2006numerical} and \cite{birgin2005numerical}.

\section{Constrained DDP using KKT conditions and slack variables}
\label{sec:CDDP}
We will henceforth focus on the constrained optimal control problem:
\begin{equation}
\label{constrained_optimalcontrol}
\begin{split}
&\min_{\bm{U}}\hspace{0.8mm}J(\bm{X},\bm{U})=\min_{\bm{U}}\big[{\textstyle\sum_{k=0}^{N-1} l(\bm{x}_k,\bm{u}_k)}+\phi(\bm{x}_N)\big]\\
&\text{subject to:}\hspace{3mm} \bm{x}_{k+1}=\bm{f}(\bm{x}_k,\bm{u}_k),\quad g_{i,k}(\bm{x}_k,\bm{u}_k)\leq0,\\
&\qquad \qquad \qquad k=0,...,N-1,\quad i=1,...,w.
\end{split}
\end{equation}
Note that we did not include equality constraints above only for compactness. Our results can be readily extended to this case as well. 

\subsection{Backward Pass}
Similar to normal unconstrained DDP, the backward pass operates on quadratic approximations of the $Q$ functions about the nominal rollouts (see eqs. \eqref{Qfunction}, \eqref{Qexpanded}, \eqref{bellman}). For the constrained case, we can write this as:
\begin{align} \label{minQ-constraints}
    &\min_{\delta\bm{u}_{k}} Q_k(\delta \bm{x}_k,\delta \bm{u}_k)\\\notag
    \text{subject to} \qquad 
    &\tilde{\bm{g}}_k(\bar{\bm{x}}_k+\delta\bm{x}_k,\bar{\bm{u}}_k+\delta\bm{u}_k) \leq 0.
\end{align}
$\tilde{\bm{g}}_k$ above is associated with the constraints influenced directly by states and controls at time instance $t_k$. We will discuss later the selection of such constraints.

We proceed by linearizing the constraints, as well as incorporating the approximate $Q$ function from \eqref{Qexpanded}. We have
\begin{align}\label{eq:g_expanded}
&\tilde{\bm{g}}(\bar{\bm{x}}_k+\delta{\bm{x}_k},\bar{\bm{u}}_k+\delta{\bm{u}_k})\\\notag
&\approx \tilde{\bm{g}}(\bar{\bm{x}}_k,\bar{\bm{u}}_k)
+ \underbrace{\tilde{\bm{g}}_{\bm u}(\bar{\bm{x}}_k,\bar{\bm{u}}_k)}
_{\bm{C}_k}\delta{\bm{u}}_k+\underbrace{\tilde{\bm{g}}_{\bm x}(\bar{\bm{x}}_k,\bar{\bm{u}}_k)}_{\bm{D}_k}\delta\bm{x}_k\notag
\end{align}
Now, for the approximate problem of \eqref{minQ-constraints}, the necessary optimality conditions from section \ref{subsec:kkt} will read as:
\begin{align}
Q_{\bm{uu}} \delta \bm{u}_k+Q_{\bm{u}} + Q_{\bm{ux}}\delta{\bm{x}_k} +\bm{C}^\mathsf{T}\bm{\lambda}_k = \bm{0}\\
{\tilde{g}_i(\bm{x}_k)\leq0, \quad \lambda}_{i,k} \geq 0, \quad \text{and} \quad \lambda_{i,k} \tilde{g}_{i} = 0 \label{kkt-eq}
\end{align}
where we have dropped the time index on the constraints and $Q$ derivatives for simplicity.

We will now rewrite the above conditions by considering a set of slack variables, such that $s_i+{g}_i=0$ and ${s}_i\geq 0$. Hence, eq. \eqref{kkt-eq} becomes 
\begin{align}
s_{i,k} \geq 0, \quad {\lambda}_{i,k} \geq 0,\quad
s_{i,k} {\lambda}_{i,k}= 0.
\end{align}
To proceed, we will consider perturbations of the slack variables and Lagrange multipliers about their (given) nominal values. Hence we obtain
\begin{align}\label{Quu-kkt}
Q_{\bm{uu}} \delta \bm{u}_k+Q_{\bm{u}} + Q_{\bm{ux}}\delta{\bm{x}_k} +\bm{C}_{k}^\mathsf{T}(\bar{\bm{\lambda}}_k+\delta\bm{\lambda}_k) = \bm{0}\\
(\bar{{s}}_{i,k} + \delta s_{i,k})(\bar{\lambda}_{i,k} + \delta{{\lambda}_{i,k}})  = {0}
\end{align}
By omitting the second-order terms, we get
\begin{align}
    \bm{S}\bar{\bm\lambda} + \bm{\Lambda}\delta\bm{s}_k + \bm{S}\delta{\bm{\lambda}_k} &= \bm{0}\\
    \text {where}\quad \bm{\Lambda}: = {\rm{diag}}(\bar{\lambda}_{i,k}), \quad \bm{S} &:= {\rm{diag}}(\bar{s}_{i,k})\notag
\end{align}
  
Moreover, the slack formulation of the inequality constraints will give
\begin{align}
    \bm{S}\bm{e} + \delta\bm{s}_k+ \tilde{\bm{g}}(\bar{\bm{x}}_k,\bar{\bm{u}}_k)
+ {\bm{C}}\delta{\bm{u}_k}+{\bm{D}}\delta\bm{x}_k = \bm{0}\\ \notag
\text{where} \qquad \bm{e} = (1,\dots, 1)
\end{align}

Overall, the obtained KKT system can be written in matrix form as
\begin{align}
    \label{eq:finalkktform}
    \begin{bmatrix}
        Q_{\bm{uu}} & \bm{0} & \bm{C}^{\mathsf{T}}\\
        \bm{0} & \bm{\Lambda} & \bm{S}\\
        \bm{C} & \bm{I} & \bm{0}
    \end{bmatrix}
    \begin{bmatrix}
    \delta \bm{u}_k\\ \delta\bm{s}_k\\ \delta\bm{\lambda}_k
    \end{bmatrix}
    &=
    \begin{bmatrix}
    -Q_{\bm{ux}}{\delta}{\bm{x}_k}-Q_{\bm{u}}-\bm{C}^\mathsf{T}\bar{\bm{\lambda}}_k\\
    -\bm{S}\bar{\bm{\lambda}}_k+ \mu_k\sigma_k\bm{e} \\
    -\bm{D}\delta\bm{x}_k-\tilde{\bm{g}}(\bar{\bm{x}}_k,\bar{\bm{u}}_k)-\bm{S}\bm{e}
    \end{bmatrix}\\\notag
    \text {with} \qquad {s}_{i,k} &\geq {0}, \quad {\lambda}_{i,k} \geq {0} .
\end{align}
We optimize this system using primal-dual interior point method \cite[chapter  18]{nocedal2006numerical}. $\bar{\bm{\lambda}}$ is initialized as $\bar{\bm{\lambda}} = \bm{e}$ ($\bm{\Lambda} = \bm{I}$) since {\it Lagrange multipliers} are required to be positive.
For slack variables $s$, they are initialized as
\begin{align}\label{eq:s_init}
 \bar{s}_{k,i}={\rm{max}}(-g_i,\epsilon),
\end{align}
where $\epsilon$ is a small positive number to keep $s_i$ positive and numerically stable. We used $\epsilon = 10^{-4}$.

In the second row of the right hand side of \eqref{eq:finalkktform}, we introduced duality measure:
\begin{align}
    \mu_k = \bar{\bm{s}}^{\mathsf{T}}_{k}\bar{\bm{\lambda}}_{k}/w.
\end{align}
It is known that if we use the pure Newton direction obtained by $\mu = 0$, we can take only a small step $\alpha$ before violating $\bm{s}^{\mathsf{T}}\bm{\lambda} \geq 0$. To make the direction less aggressive, and the optimization process more effective we reduce $s_{i}\lambda_{i}$ to a certain value based on the average value of elementwise product $s_{i}\lambda_{i}$, instead of zero. Note that $\mu$ is an average value of $s_{i}\lambda_{i}$ and $\mu$ must converge to zero over the optimization process. We satisfy this requirement by multiplying $\sigma$ ($0<\sigma<1$) \cite[Chapter 19]{nocedal2006numerical}.\\
$\sigma$ is given by:
\begin{align}
    \sigma_k &= 0.1 \min(0.05\frac{1-\xi_k}{\xi_k},2)^3\\
    {\text {where}} \quad \xi_k &= \frac{\min_i (s_{i,k} \lambda_{i,k})}{\mu_k}.
\end{align}

Our goal is to solve the above system analytically, which, as we will find, requires the inversion of $\bm{S}$ and $\bm{\Lambda}$. It might be the case, however, that these matrices are close to being singular; for example, elements of $\bm{S}$ will be close to zero when the corresponding constraints approach their boundaries. To tackle this problem we will perform the following change of variables
\begin{align}
    \delta\bm{p}_k := \bm{S}^{-1} \delta \bm{s}_k, \qquad
    \delta \bm{q}_k := \bm{\Lambda}^{-1} \delta \bm{\lambda}_k.
\end{align}
Then the new KKT system can be obtained as
\begin{align}
    \nonumber
    \begin{bmatrix}
        Q_{\bm{uu}} & \bm{0} & \bm{C}^{\mathsf{T}}\bm{\Lambda}\\
        \bm{0} & \bm{\Lambda}\bm{S} & \bm{S}\bm{\Lambda}\\
        \bm{\Lambda}\bm{C} & \bm{\Lambda}\bm{S} & \bm{0}
    \end{bmatrix}
    \begin{bmatrix}
    \delta \bm{u}_k\\ \delta\bm{p}_k\\ \delta\bm{q}_k
    \end{bmatrix}\qquad\qquad\qquad\qquad\qquad\\
\qquad =
    \begin{bmatrix}
    -Q_{\bm{ux}}{\delta}{\bm{x}_k}-Q_{\bm{u}}-\bm{C}_k^\mathsf{T}\bar{\bm{\lambda}}_k\\
    -\bm{S}\bar{\bm{\lambda}}_k + \mu_k\sigma_k\bm{e}\\
    -\bm{\Lambda}(\bm{D}_k\delta\bm{x}_k+\tilde{\bm{g}}(\bar{\bm{x}}_k,\bar{\bm{u}}_k)+\bm{S}\bm{e})
    \end{bmatrix}
    =
    \begin{bmatrix}
    \bm{a}\\
    \bm{b}\\
    \bm{\Lambda d}\\
    \end{bmatrix}\\
    \notag
    \text {s.t.} \qquad s_{i,k} \geq 0, \quad \lambda_{i,k} \geq 0. 
\end{align}
We can now avoid singularity issues, since $\bm{\Lambda S} \rightarrow \mu_k \sigma_k$ (instead of 0) with our new formulation.

Now notice that in the backward pass of DDP, we do not have $\delta\bm{x}$. Hence, our strategy will be to first solve the KKT system on the nominal rollout by substituting $\delta\bm{x} = \bm{0}$, and then use our optimal values of $\bm s$ and $\bm \lambda$ as our $\bar{\bm s}$ and $\bar{\bm \lambda}$ for the next KKT iteration.\\
An analytical solution for this case can be obtained by algebraic manipulations, which we state below:
\begin{align}
\delta\bm{q} &= \bm{M}(\bm{E}\bm{Q}_{uu}^{-1}\bm{a}-\bm{\Lambda}\bm{d}+\bm{b})\\
\delta\bm{p} &= \bm{F}^{-1}(\bm{b}-\bm{F}\delta\bm{q})\\
\delta \bm{u} &= \bm{a}-\bm{E}^{\mathsf{T}}\delta\bm{q}\\
\text {where} \qquad \bm{M} &= (\bm{E}\bm{Q}_{uu}^{-1}\bm{E}^{\mathsf{T}}+\bm{F})^{-1}\\\notag
\bm{E}&= \bm{\Lambda}\bm{C}, \quad \bm{F} = \bm{\Lambda}\bm{S}.\label{eq:EF}
\end{align} 
We will next update $\bar{\bm{s}}$ and $\bar{\bm{\lambda}}$ by
\begin{align}
\bar{\bm{s}}_k = \bar{\bm{s}}_k + \alpha \delta {\bm{s}}_k, \quad
\bar{\bm{\lambda}}_k = \bar{\bm{\lambda}}_k + \alpha \delta{\bm{\lambda}_k}.
\end{align}
The step size $\alpha$ must be determined to keep $\bm{s}_k$ and $\bm{\lambda}_k$ non-negative. The following strategy is also used in interior point methods \cite[Chapter 16]{nocedal2006numerical}.
\begin{align}
\alpha = \min(\alpha_{s}, \alpha_\lambda)\\
\alpha_s = \min_{\delta s_i^{k} < 0}(1,-\zeta\frac{s_{i,k}}{\delta s_{i,k}}), \quad
\alpha_\lambda &= \min_{\delta \lambda_{i,k} < 0}(1,-\zeta\frac{\lambda_{i,k}}{\delta \lambda_{i,k}})\\
{\text{where}}\quad 0.9 &\leq \zeta < 1. \notag
\end{align}
When there are no negative elements in $\delta\bm{s}$ or $\delta\bm{\lambda}$, the corresponding step sizes are taken to be $\alpha_{s} = 1$ and $\alpha_{\lambda} = 1$.
Using this step size, we also update the linearized constraint function and the $Q$ function. For convenience, we write the new $\bar{\bm{u}}_k$ as $\bar{\bm{u}}_k + \alpha\delta{\bm{u}}_k^{j-1}$. $j-1$ implies that $\delta \bm{u}$ is form one iteration before. Then the updated constraint function on the nominal trajectory $(\delta \bm{x} = \bm{0})$ is:
\begin{align}
    &\quad\tilde{\bm{g}}(\bar{\bm{u}}_k + \alpha\delta{\bm{u}}_k^{j-1} + \delta{\bm{u}_k})\\\notag
    &\approx
    \tilde{\bm{g}}(\bar{\bm{u}}_k) + \tilde{\bm{g}}_{u}(\alpha \delta{\bm{u}}_k^{j-1} + \delta{\bm{u}_k})\\\notag
    &=(\tilde{\bm{g}}(\bar{\bm{u}}_k) + \bm{C}_k\alpha \delta{\bm{u}}_k^{j-1} ) +  \bm{C}_k\delta{\bm{u}_k}.
    \end{align}
    Therefore the updated nominal constraint function is
    \begin{align}
        \bm{g}(\bar{\bm{u}}_k) = \bm{g}(\bar{\bm{u}}_k) + \bm{C}_k\alpha \delta{\bm{u}_k^{j-1}}.
    \end{align}
    For $Q$ function, we expand them around nominal trajectory considering small perturbation.
    \begin{align}
    &{Q}(\bar{\bm{u}}_k + \alpha\delta{\bm{u}}_k^{j-1} + \delta\bm{u}_k)\\ \notag &\approx
    \frac{1}{2}(\alpha \delta{\bm{u}_k^{j-1}} + \delta{\bm{u}_k})^{\mathsf{T}}{Q_{\bm{uu}}}(\alpha \delta{\bm{u}_k^{j-1}} + \delta{\bm{u}_k})\\\notag
    &\qquad+ Q_{\bm{u}}^{\mathsf{T}}(\alpha \delta{\bm{u}_k^{j-1}} + \delta{\bm{u}_k}).
\end{align}
Using this new $\bm{Q}$ function we construct {\it Lagrangian} as,
\begin{align}
    &L = \frac{1}{2}(\alpha \delta{\bm{u}_k^{j-1}} + \delta{\bm{u}_k})^{\mathsf{T}}{Q_{\bm{uu}}}(\alpha \delta{\bm{u}_k^{j-1}} + \delta{\bm{u}_k})\\\notag & +Q_{\bm{u}}^{\mathsf{T}}(\alpha \delta{\bm{u}_k^{j-1}} + \delta{\bm{u}_k}) + \bm{\lambda}^{\mathsf{T}}(\tilde{\bm{g}}(\bar{\bm{u}}_k) + \bm{C}_k\alpha \delta{\bm{u}}_k^{j-1}  +  \bm{C}_k\delta{\bm{u}_k}).
\end{align}
From KKT condition $\frac{\partial L}{\partial \delta\bm{u}_k}$, we have
\begin{align}
    Q_{\bm{uu}} (\delta \bm{u}_k + \alpha \delta \bm{u}^{j-1})+Q_{\bm{u}}+\bm{C}_{k}^\mathsf{T}\bm{\lambda} = \bm{0}.
\end{align}
Comparing the above equation with \eqref{Quu-kkt} on the nominal trajectory ($\delta{\bm{x}_k} = \bm{0}$), the new $Q_{\bm{u}}$ can be obtained as
\begin{align}
    Q_{\bm{u}} = Q_{\bm{u}} + \alpha Q_{\bm{uu}} \delta \bm{u}^{j-1}.
\end{align}
and the KKT system \eqref{eq:finalkktform} is iteratively solved, updating $\bar{\bm{s}}, \bar{\bm{\lambda}}$, and $\mu_k$, until the duality measure $\mu_k$ is improved to a certain threshold. In this paper we used 0.01 as the threshold. \\
Finally, using the updated $\bar{\bm{s}}_k$ and $\bar{\bm{\lambda}}_k$, we solve the system at the perturbed trajectory $(\delta \bm{x} \neq 0)$ for $\delta\bm{u}_k$:
\begin{align}
    {\delta}{\bm{u}_k} &= -{Q_{\bm{uu}}^{-1}}[\bm{H}Q_u+\bm{E}^\mathsf{T}{\bm{M}}{\bm{\Lambda}}(\tilde{\bm{g}}(\bar{\bm{x}}_k, \bar{\bm{u}}_k)+\bm{Se})]\nonumber \\
    &-{Q_{\bm{uu}}^{-1}}({\bm{H}}{Q_{\bm{ux}}} + {\bm{E}}^{\mathsf{T}}{\bm{M}}\bm{\Lambda}\bm{D}){\delta}{\bm{x}_k}\\
    \text{where} \qquad {\bm{H}} &= \bm{I}-\bm{E}^{\mathsf{T}}\bm{M}\bm{E}{Q_{\bm{uu}}^{-1}}\label{eq:H}.
\end{align}
The feedforward gain, $\bm{k}$, and feedback gain, ${\bm{K}}$, can then be obtained as follows:
\begin{align}
    \bm{k}&= -{Q_{\bm{uu}}^{-1}}[\bm{H}\bm{Q}_u+\bm{E}^\mathsf{T}{\bm{M}}{\bm{\Lambda}}(\bm{g}(\bar{\bm{x}}_k, \bar{\bm{u}}_k)+\bm{Se})]\label{eq:k_constrained}\\
    \bm{K} &= -{Q_{\bm{uu}}^{-1}}({\bm{H}}{Q_{\bm{ux}}} + {\bm{E}}^{\mathsf{T}}{\bm{M}}\bm{\Lambda}\bm{D})\label{eq:K_constrained}
\end{align}
We will finally discuss which constraints $\tilde{\bm g}$ to consider. We assume the full state $\bm{x}$ is composed of position $\bm{x}^p {\in}{\mathbb{R}}^{n_p}$, a function of state itself, and $\bm{x}^v {\in}{\mathbb{R}}^{n_v}$, a function of state and control: $\bm{x} = [\bm{x}^{p\mathsf{T}},\bm{x}^{v\mathsf{T}}]^{\mathsf{T}}$. p implies position, and v implies velocity.
\begin{align}
    \bm{x}_{k+2}^p = \bm{f}^p(x_{k+1}), \quad
    \begin{bmatrix}
        \bm{x}_{k+1}^p\\
        \bm{x}_{k+1}^v
    \end{bmatrix}
    =\begin{bmatrix}
    \bm{f}^p(\bm{x}_k)\\
    \bm{f}^v(\bm{x}_k, \bm{u}_k)
    \end{bmatrix}.
\end{align}
$\bm{u}_k$ (control at time step $k$) that we obtain by solving QP does not affect the $\bm{x}_k$ state at the same time step, but $\bm{x}^p$ two time steps forward $\bm{x}_{k+2}^p$, and $\bm{x}^v$ one time step forward $\bm{x}^v_{k+1}$. This implies that we should solve the QP subject to constraints two time steps forward for $\bm{x}^p$, and one time step forward for $\bm{x}^v$.
We use this fact to resolve the problem when $\bm{C}$ is zero in our algorithm. If $\bm{C}$ is zero, feedback and feedfoward gains are the same as normal unconstrained DDP, see \eqref{eq:EF}, \eqref{eq:H},  \eqref{eq:k_constrained}, and \eqref{eq:K_constrained}.
First, we divide elements of $\bm{g}_k$ into function of $\bm{x}^p$, $\bm{x}^v$, and $\bm{u}_k$, and write them as $\bm{g}^p_{k}$, $\bm{g}^v_{k}$, and $\bm{g}^c_{k}$.
For $\bm{g}^{p}_k$ we propagate two time steps, and for $\bm{g}^v_k$ one time step forward to make them explicit functions of $\bm{u}_k$.   
Expanded $\bm{g}^p_{k+2}$ can be calculated using the chain rule:
\begin{align}
&\bm{g}^p_{k+2}(\bar{\bm{x}}_k+\delta{\bm{x}_k},\bar{\bm{u}}_k+\delta{\bm{u}_k})\\\notag
&\approx\bm{g}^{p}_{k+2}(\bar{\bm{x}}_k,\bar{\bm{u}}_k)+
\underbrace{\frac{\partial \bm{g}^{p}_{k+2}}{\partial {\bm{x}}_k}}_{\bm{D}^p_k}
\delta {\bm{x}}_{k} +
\underbrace{\frac{\partial \bm{g}^{p}_{k+2}}{\partial{\bm{u}}_k}}_{\bm{C}^p_k}{\delta\bm{u}_k} 
\end{align}
where
\begin{align}
    \frac{\partial {\bm{g}}^p_{k+2}}{\partial{\bm{x}}_k} =
    \frac{\partial {\bm{g}}_{k+2}^{p}}{\partial {\bm{x}}^p_{k+2}}
    \frac{\partial {\bm{x}}^p_{k+2}}{\partial {\bm{x}}_{k+1}}
    \frac{\partial {\bm{x}}_{k+1}}{\partial {\bm{x}}_{k}}
    =\frac{\partial {\bm{g}}_{k+2}^{p}}{\partial {\bm{x}}^p_{k+2}}
    \bm{f}^{p}_{\bm{x},k+1}
    \bm{f}_{\bm{x},k}\\
    \frac{\partial {\bm{g}}^p_{k+2}}{\partial {\bm{u}}_k} =
    \frac{\partial {\bm{g}}_{k+2}^{p}}{\partial {\bm{x}}_{k+2}^p}
    \frac{\partial {\bm{x}}^p_{k+2}}{\partial {\bm{x}}_{k+1}}
    \frac{\partial {\bm{x}}_{k+1}}{\partial {\bm{u}}_{k}}
    = \frac{\partial {\bm{g}}_{k+2}}{\partial {\bm{x}}^p_{k+2}}
    \bm{f}^{p}_{\bm{x},k+1}
    \bm{f}_{\bm{u},k}
\end{align}
$\bm{g}^v_{k+1}$ is:
\begin{align}
&\bm{g}^v_{k+1}(\bar{\bm{x}}_k+\delta{\bm{x}_k},\bar{\bm{u}}_k+\delta{\bm{u}_k})\\\notag
&\approx\bm{g}^{v}_{k+1}(\bar{\bm{x}}_k,\bar{\bm{u}}_k)+
\underbrace{\frac{\partial \bm{g}^{v}_{k+1}}{\partial {\bm{x}}_k}}_{\bm{D}^v_{k}}
\delta {\bm{x}}_{k} +
\underbrace{\frac{\partial \bm{g}^{v}_{k+1}}{\partial{\bm{u}}_k}}_{\bm{C}^v_k}{\delta\bm{u}_k} 
\end{align}
\begin{align}
    \frac{\partial {\bm{g}}^{v}_{k+1}}{\partial{\bm{x}}_k} =
    \frac{\partial {\bm{g}}_{k+1}^{v}}{\partial {\bm{x}}_{k+1}}
    \frac{\partial {\bm{x}}_{k+1}}{\partial {\bm{x}}_{k}}
    =\frac{\partial {\bm{g}}_{k+1}^{v}}{\partial {\bm{x}}_{k+1}}
    \bm{f}_{\bm{x},k}\\
    \frac{\partial {\bm{g}}^v_{k+1}}{\partial {\bm{u}}_k} =
    \frac{\partial {\bm{g}}_{k+1}^{v}}{\partial {\bm{x}}_{k+1}}
    \frac{\partial {\bm{x}}_{k+1}}{\partial {\bm{u}}_{k}}
    = \frac{\partial {\bm{g}}_{k+1}^{v}}{\partial {\bm{x}}_{k+1}}
    \bm{f}_{\bm{u},k}
\end{align}
For $\bm{g}^c$ we use the same expressions as $\tilde{\bm{g}}$.

Stacking $\bm{g}$, $\bm{C}$, and $\bm{D}$, we have the linearized constraints
\begin{align}
    \bm{g}(\bar{\bm{x}}_k,\bar{\bm{u}}_k)
    = \begin{bmatrix}
    \bm{g}^p_{k+2}\\\bm{g}^v_{k+1}\\\bm{g}^c_{k}
    \end{bmatrix},
    \bm{C}_k = \begin{bmatrix}
    \bm{C}^p_k\\\bm{C}^v_k\\\bm{C}^c_k\\
    \end{bmatrix},
    \bm{D}_k = \begin{bmatrix}
    \bm{D}^p_k\\\bm{D}^v_k\\\bm{D}^c_k\\
    \end{bmatrix}
\end{align}
 \begin{algorithm}
 \caption{Backward Pass}
 \begin{algorithmic}[1]
 \label{alg:backpass}
 \renewcommand{\algorithmicrequire}{\textbf{Input:}}
 \renewcommand{\algorithmicensure}{\textbf{Output:}}
  \STATE {Initialize:
  $V_{N} \gets \phi(\bm{\bar{x}}_{N})$\\
  $V_{\bm{x},N} \gets \bm{\nabla_x}\phi(\bm{\bar{x}}_{N})$,
  $V_{\bm{xx},N} \gets \bm{\nabla_{xx}}\phi(\bm{\bar{x}}_{N})$}
  \FOR {$k = N-1$ to $0$}
  \STATE {Calculate $l,Q$, and their derivatives at $k$}
  \STATE {Regularize $Q_{\bm{uu}}$,$Q_{\bm{ux}}$, and $Q_{\bm{ux}}$ using $\nu_1$ and $\nu_2$}
\\ \textit{$\bar{\bm{s}}$ and $\bar{\bm{\lambda}}$ process}
  \STATE {Initialize $\bar{s}_{k,i}={\rm{max}}(-g_i,\epsilon)$,
  $\bar{\lambda}_{k,i} =1$}
  \STATE {$\mu_0 \gets \bm{s}^{\mathsf{T}}\bm{\lambda}/w$}
  \WHILE{$\mu/\mu_0 < 0.01$}
  \STATE{Solve KKT system to obtain $\delta\bm{u}_k$, $\delta\bm{s}_k$, and $\delta\bm{\lambda}_k$}
  \STATE{Determine step size $\alpha$}
  \STATE{Update $\bar{\bm{u}}_k, \bar{\bm{s}}_k, \bar{\bm{\lambda}}_k, \bm{g}(\bar{\bm{x}}_k, \bar{\bm{u}}_k)$, and $Q_{\bm{u}}$ }
  \STATE {$\mu \gets \bar{\bm{s}}_k^{\mathsf{T}}\bar{\bm{\lambda}}_k/w$}
  \ENDWHILE
  \\ \textit{Update gains and value functions}
  \STATE {$\bm{k} \gets  -Q_{\bm{uu}}^{-1}[\bm{H}\bm{Q}_u+\bm{E}^\mathsf{T}{\bm{M}}{\bm{\Lambda}}(\bm{g}(\bar{\bm{x}}, \bar{\bm{u}})+\bm{Se}]$}\\
  \STATE {$\bm{K} \gets -Q_{\bm{uu}}^{-1}({\bm{H}}{Q_{\bm{ux}}} +{\bm{E}}^{\mathsf{T}}{\bm{M}}\bm{\Lambda}\bm{D})$}
  \STATE  {$V_{\bm{x}} \gets Q_{\bm{x}} + \bm{K}^{\mathsf{T}}Q_{\bm{uu}}\bm{k}
    +\bm{K}^{\mathsf{T}}Q_{\bm{u}} + Q_{\bm{ux}}^{\mathsf{T}}\bm{k}$}
    \STATE {$V_{\bm{xx}} \gets Q_{\bm{xx}} + \bm{K}^{\mathsf{T}}Q_{\bm{uu}}\bm{K}
    +\bm{K}^{\mathsf{T}}Q_{\bm{ux}} + Q_{\bm{ux}}^{\mathsf{T}}\bm{K}$}
  \ENDFOR
 \STATE {Store Derivatives of $Q$} 
 \end{algorithmic} 
 \end{algorithm}

 \begin{algorithm}
 \caption{Algorithm for forward pass}
 \begin{algorithmic}[1]
 \renewcommand{\algorithmicrequire}{\textbf{Input:}}
 \renewcommand{\algorithmicensure}{\textbf{Output:}}
 \newcommand{\algorithmicbreak}{\textbf{break}}
\newcommand{\BREAK}{\STATE \algorithmicbreak}
  \STATE {Calculate $J_{int} \gets J(\bm{X},\bm{U})$}
  \STATE {$\bm{x} \gets \bm{x}_0$}
  \STATE {$\bm{X}_{\rm temp}\gets\bm{X},\quad \bm{U}_{\rm temp}\gets\bm{U}$}
    \FOR {$k = 0$ to $N-2$}
  \STATE {$\delta \bm{x} \gets \bm{x}-\bm{x}_k$}
  \STATE {$\bm{x}_{{\rm{temp}},k} \gets \bm{x}$}
  \\ \textit{Solve QP:}
  \STATE {${\delta \bm{u}_{k}^{\ast}} = \argmin[\frac{1}{2}\delta\bm{u}_{k}^{\mathsf{T}}Q_{\bm{uu}}\delta\bm{u}_k+
    \delta\bm{u}_{k}^{\mathsf{T}}(Q_{\bm{u}}+Q_{\bm{ux}}\delta\bm{x}_k)]$, subject to ${\bm{g}_k}({\bm{x}},\bm{u}_k)+ {\bm{C}_{k}}({\bm{x}},\bm{u}_k){\delta}{\bm{u}_k}\leq \bm{0}$}
    \STATE{$\bm{u}_{\rm{temp},k} = \bm{u}_{k}+\delta \bm{u}_{k}^{\ast}$}
    \STATE{$\bm{x}\gets\bm{f}(\bm{x},\bm{u}_{\rm{temp},k})$}
    \ENDFOR
    \STATE {$\delta \bm{x} \gets \bm{x}-\bm{x}_{N-1}$}
  \STATE {$\bm{x}_{{\rm{temp}},N-1} \gets \bm{x}$}
  \WHILE{$\rm{flag = \textbf{False}}$}
  \STATE{$\rm{flag = \textbf{True}}$}
  \\ \textit{Solve QP:}
  \STATE {${\delta \bm{u}^{\ast}_{N-1}} = \argmin[\frac{1}{2}\delta {\bm{u}_{N-1}}^{\mathsf{T}}Q_{\bm{uu}}\delta\bm{u}_{N-1}+\delta\bm{u}_{N-1}^{\mathsf{T}}(Q_{\bm{u}}+Q_{\bm{ux}}\delta\bm{x}_{N-1})]$, subject to $|\bm{u}_{N-1}| \leq \varDelta$}
     \\ \textit{Update $\bm{x}, \bm{u}$ and Check Feasibility of QP}
    \STATE{$\bm{u}_{{\rm temp},N-1} = \bm{u}_{N-1}+\delta \bm{u}_{N-1}^{\ast}$}
    \STATE{$\bm{x}_{{\rm temp},N}\gets\bm{f}(\bm{x},\bm{u}_{\rm{temp},N-1})$}
    \IF{${\textbf{any}}\quad {\bm{g}}(\bm{x}_{{\rm temp},N})>0$}
    \STATE{$\rm{flag = \textbf{False}}$}
    \STATE{$\varDelta \gets \eta\varDelta$}
    \BREAK
    \ENDIF
    \ENDWHILE  
  \STATE {Calculate $J_{temp} \gets J(\bm{X}_{\rm{temp}},\bm{U}_{\rm{temp}})$}

  \IF{$J_{\rm{temp}} < J_{\rm{int}}$}
  \STATE {$\bm{X}\gets\bm{X}_{\rm{temp}}, \bm{U}\gets\bm{U}_{\rm{temp}}$}
  \STATE {Decrease $\nu_1$ and $\nu_2.$}
  \ELSE
  \STATE {Increase $\nu_1$ and $\nu_2.$}
  \ENDIF
 \end{algorithmic} 
 \end{algorithm}

\subsection{Forward Pass}
The following QP problem is solved which takes into account all the constraints to guarantee the updated nominal trajectory is feasible:
\begin{align}{\label{forawrd_QP}}
    \argmin_{\delta \bm{u}_k}[\frac{1}{2}\delta {\bm{u}_k}^{\mathsf{T}}&Q_{\bm{uu}}\delta\bm{u}_k
    +
    Q_{\bm{u}}^{\mathsf{T}}\delta\bm{u}_{k}+
    \delta {\bm{u}_k}^{\mathsf{T}}Q_{\bm{ux}}\delta\bm{x}_k]\\\notag
\text{subject to} \quad 
        &\quad{\bm{g}_k}(\bar{\bm{x}}_{k}+\delta{\bm{x}_{k}},\bar{\bm{u}}_{k}+\delta{\bm{u}_{k}}) \quad\\\notag
        &\approx {\bm{g}_k}(\bar{\bm{x}}_{k},\bar{\bm{u}}_{k})+\bm{D}_{k}{\delta}{\bm{x}}_{k}+ {\bm{C}_{k}}{\delta}{\bm{u}_{k}} \leq \bm{0}
\end{align}
Note that $\delta{\bm{u}}_k$ is used to obtain $\bm{x}_{k+1}$ and we already have $\delta{\bm{x}}_k$. The updated nominal state, $\bar{\bm{x}}'_k = \bar{\bm{x}}_k + \delta \bm{x}_k$, can be used instead of the previous iteration's $\bar{\bm{x}}_k$ and $\delta \bm{x}_k$. Thus the $\delta \bm{x}$ term for $\bar{\bm{x}}'_k$ in the linearized constraint equation becomes zero and we obtain the following:
\begin{align}\label{eq:linearized_constraint_forward}
    {\bm{g}_k}(\bar{\bm{x}}_k',\bm{u}_k+\delta{\bm{u}}_k) = {\bm{g}_k}(\bar{\bm{x}}_k',\bm{u}_k)+ {\bm{C}_{k}}(\bar{\bm{x}}_k',\bm{u}_k){\delta}{\bm{u}_k}\leq \bm{0}
\end{align}
Again, one or two time step propagation of $\bm{g}$ shown in backward pass is important, because otherwise $\bm{C}_k$ might be zero, and $\bm{u}_k$ does not show up in the linearized constraint. At time step $N$, constraints $\bm{g}^{p}_{N+2}$ are not available because the time horizon ends at time step $N+1$. Hence, we solve the QP under box constraints
\begin{align}
    -\bm{\varDelta} \leq \delta \bm{u}_N \bm{\leq \varDelta}
\end{align}
$\bm{\varDelta}$ is a vector of an adaptive trust region initialized by a relatively large positive value. If the solution is not feasible, $\bm{\varDelta}$ is made to be smaller by
\begin{align}\label{eq:trust_region}
 \bm{\varDelta} = \eta \bm{\varDelta}, \qquad \text{where}\quad 0 < \eta < 1   
\end{align}
   
This makes the new trajectory closer to the trajectory one iteration before, which is feasible. The trust region is made smaller repeatedly until the solution of the QP becomes feasible.

\subsection{Regularization}

In DDP, regularization plays a big role and highly affects the convergence.
We use the regularization scheme and scheduling technique proposed in \cite{TassaDDP2012}.
\begin{align}
Q_{\bm{xx}} &=\bm{l_{xx}}+{\bm{f}}_{\bm{x}}^{\mathsf{T}}(V_{\bm{xx}}' + \nu_1 \bm{I}_n){\bm{f}}_{\bm{x}}\\
 Q_{\bm{ux}} 
     &=\bm{l_{ux}}+\bm{f}_{\bm{u}}^{\mathsf{T}}(V_{\bm{xx}}' + \nu_1 \bm{I}_n)\bm{f}_{\bm{x}}\\
    Q_{\bm{uu}} &= \bm{l_{uu}}+\bm{f}_{\bm{u}}^{\mathsf{T}}(V_{\bm{xx}}'+ \nu_1 \bm{I}_n)\bm{B} + \nu_2 \bm{I}_m
\end{align}
These regularized derivatives are used to calculate gains $\bm{k}$ and $\bm{K}$, and instead of using \eqref{eq:value_riccati} we update the value function as follows:
\begin{align}
    V_{\bm{x}} &= Q_{\bm{x}} + \bm{K}^{\mathsf{T}}Q_{\bm{uu}}\bm{k}
    +\bm{K}^{\mathsf{T}}Q_{\bm{u}} + Q_{\bm{ux}}^{\mathsf{T}}\bm{k}\\
    V_{\bm{xx}} &= Q_{\bm{xx}} + \bm{K}^{\mathsf{T}}Q_{\bm{uu}}\bm{K}
    +\bm{K}^{\mathsf{T}}Q_{\bm{ux}} + Q_{\bm{ux}}^{\mathsf{T}}\bm{K}
\end{align}
Both $\nu_1$ and $\nu_2$ are positive values. $\nu_1$ makes the new trajectory closer to the previous one, while $\nu_2$ makes the control step conservative and makes ${Q}_{\bm{uu}}$ positive definite.
\section{Constrained DDP with Augmented Lagrangian methods} \label{sec:AL-DDP}
\subsection{AL-DDP and penalty function}
Here, we will be using the Augmented Lagrangian approach to extend DDP for solving \eqref{constrained_optimalcontrol}. We call this technique AL-DDP. The main idea is to observe that {\it partial elimination of constraints} can be used on the inequality constraints $\bm g$ of \eqref{constrained_optimalcontrol}. This means that the penalty function $\mathcal{P}$ from section \ref{subsec:AL}
can only be applied to the inequality state constraints, while the dynamics are implicitly satisfied due to DDP's problem formulation. We will thus be considering the following problem
\begin{equation}
\label{eq:constrained_control_multipliers}
\begin{split}
&\min_{\bm U}\hspace{1mm} \bigg[\underbrace{J(\bm X, \bm U) + \sum_{i,k}\mathcal{P}(\lambda_i^k, \mu_i^k, g_{i,k}(\bm{x}_k,\bm{u}_k))}_{L_A}\bigg]\\
\text{s.t.}&\hspace{7mm} \bm{x}_{k+1}=\bm{f}(\bm{x}_k, \bm{u}_k),\quad x^0=\bar{x}^0,\quad k=0, 1, ..., H-1,\\
\end{split}
\end{equation}
where $\lambda_i^k$, $\mu_i^k$ denote Lagrange multipliers and penalty parameters respectively. We will thus be using the approach discussed in section \ref{subsec:AL}, using specifically unconstrained DDP to optimize \eqref{eq:constrained_control_multipliers}, followed by an update on the Lagrange multipliers and the penalty parameters.

Since DDP requires $L_A(\cdot)$ to be twice differentiable, we selected the penalty function as
\[\mathcal{P}(\lambda_{i,k}, \mu_{i,k}, g_{i,k}(\bm{x}_k))=\frac{(\lambda_{i,k})^2}{\mu_{i,k}}\phi\bigg({\frac{\mu_{i,k}}{\lambda_{i,k}}g_{i,k}(\bm{x}_k)}\bigg),\]
with
\[\phi(t):=\begin{cases}
\frac{1}{2}t^2 + t,\quad t\geq-\frac{1}{2}\\
-\frac{1}{4}\log(-{2t})-\frac{3}{8}, \quad \text{otherwise},
\end{cases}\]
which can be viewed as a smooth approximation to the Powell-Hestenes-Rockafellar method \cite{birgin2005numerical}.
\subsection{Combination of AL-DDP and KKT frameworks}\label{sec:ALandKKT}
The Augmented Lagrangian approach is typically robust to initializations of the algorithm, but displays oscillatory behavior near the (local) solution \cite{birgin2005numerical}. This can be readily explained from optimization theory, since Multiplier methods generally converge only linearly \cite{birgin2005numerical}.

The idea here is to combine the two approaches: We begin by using the AL-DDP formulation until a pre-specified precision of the cost and constraints, and subsequently switch to the KKT-based approach of section \ref{sec:CDDP}. If sufficient improvement is not observed within a few iterations, we switch back to the Augmented Lagrangian method, and reduce the aforementioned tolerances for the ``switching" mechanism.
We also have one more reason for the combination. 
By applying the control limited DDP technique \cite{tassa2014control} in the backward pass of the Augmented Lagrangian method, we can handle the control limits as a hard box constraint. In the control limited DDP method, the feedforward gain $\bm{k}_k$ is obtained as
\begin{align}
\bm{k}_k = \argmin_{\delta \bm{u}_k}\frac{1}{2}\delta\bm{u}_{k}^{\mathsf{T}}Q_{\bm{uu}}\delta\bm{u}_{k} + \delta\bm{u}_{k}^{\mathsf{T}}Q_{\bm{u}}\\\notag
\bm{u}_{l} \leq \bm{u}_k + \delta \bm{u}_k \leq  \bm{u}_{u},
\end{align}
where $\bm{u}_{l}$ is the lower and $\bm{u}_{u}$ is the upper limit of control. For the feedback gain $\bm{K}_k$, corresponding rows to the active control limits are set to be zero.
 As we will discuss in Section \ref{sec:CDDP}, our KKT-based method can not handle a situation when state and control constraints conflict with each other. This situation typically happens when initial state is far from desired state and a large control is required. By providing a good initial trajectory from AL for the KKT-based method, our method can successfully handle both state and control constraints.
\section{Results} \label{sec:results}
 In this section we provide simulation results and comparisons between our method and prior work. We call our method ``S-KKT" named after  the slack variable and KKT conditions. We also test a slightly different version from S-KKT, in which we use the active set method \cite{Lin91} instead of slack variables, but still use one and two time step forward expansion of the constraint function. More precisely, in this method, constraints are took into account only when they are close to active, and they are linerarized under an assumption that active constraints in current iteration remains to be active in next iteration, that is
 \begin{align}\label{eq:const_assumption}
 {\bm{g}}(\bar{\bm{x}}_k+\delta{\bm{x}_k},\bar{\bm{u}}_k+\delta{\bm{u}_k})
= {\bm{g}}(\bar{\bm{x}}_k,\bar{\bm{u}}_k) \approx \bm{0}.
\end{align}
Using this assumption in \eqref{eq:g_expanded}, constraints are written as
\begin{align}
{\bm{C}_k}\delta{\bm{u}}_k+{\bm{D}_k}\delta\bm{x}_k = \bm{0},
\end{align}
and QP is solved under this equality constraints instead of \eqref{eq:linearized_constraint_forward}.
 The purpose of showing this algorithm here is to see the effect of slack variables and the assumption in \eqref{eq:const_assumption}. We call this method the ``active set method".
 We evaluate these methods and compare them with the former method \cite{xie2017differential}. The next step is to combined our methods with other optimization algorithms and evaluate the performance. 
\subsection{S-KKT}
We evaluate our constrained DDP algorithms in two different systems, a simplified 2D car and a quadrotor.
\subsubsection{2D car}\label{subsubsec:2Dcar}
We consider a 2D car following with dynamics give by the expression below \cite{xie2017differential}:
\begin{align}
       \begin{bmatrix}
     x_{k+1}\\
     y_{k+1}\\
     \theta_{k+1}\\
     v_{k+1}\\
     \end{bmatrix}
     =
     \begin{bmatrix}
     x_{k}+v_{k}\sin{\theta_k} \varDelta t\\
     y_{k}+v_{k}\cos{\theta_k} \varDelta  t\\
     \theta_{k}+u_{k}^{\theta}v_{k} \varDelta  t\\
     v_{k}+u_k^v \varDelta  t\\
    \end{bmatrix}  
\end{align}.
The car has state $\bm{x} = [x,y,\theta,v]^{\mathsf{T}}$, and control $u^{\theta}$ on the steering angle and $u_v$
on the velocity. We consider a reaching task to $\bm{x}_g = [3,3,\frac{\pi}{2},0]^{\mathsf{T}}$ while avoiding three obstacles.The obstacles are formulated as
\begin{align}
\bm{g}_{\rm{car}} =
\begin{bmatrix}
0.5^2-(x-1)^2 - (y-1)^2\\
0.5^2-(x-1)^2 - (y-2.5)^2\\
0.5^2-(x-2.5)^2 - (y-2.5)^2
\end{bmatrix} \leq \bm{0}.
\end{align}

\begin{figure}[h]
 \begin{minipage}[b]{0.34\linewidth}
  \centering
  \includegraphics[width=\columnwidth]{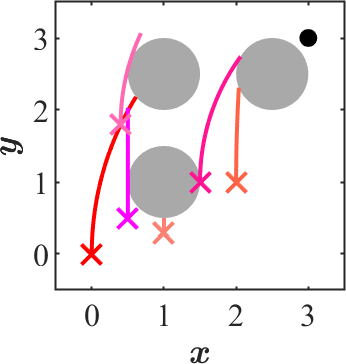}
  \subcaption{Former method}
 \end{minipage}
 \begin{minipage}[b]{0.31\linewidth}
  \centering
   \includegraphics[width=\columnwidth]{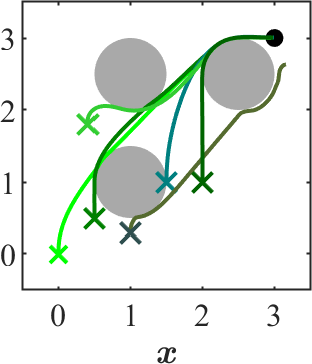} \subcaption{Active set }
 \end{minipage}
 \begin{minipage}[b]{0.31\linewidth}
  \centering
   \includegraphics[width=\columnwidth]{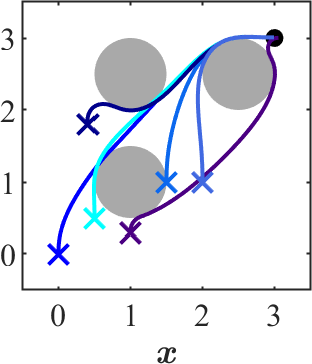} \subcaption{S-KKT}
 \end{minipage}
 \caption{2D car trajectories starting from several initial points. Starting points are shown by ``x".}\label{fig:2dcar_init_pnt}
\end{figure}

\begin{figure}[h]
 \begin{minipage}[b]{0.34\linewidth}
  \centering
  \includegraphics[width=\columnwidth]{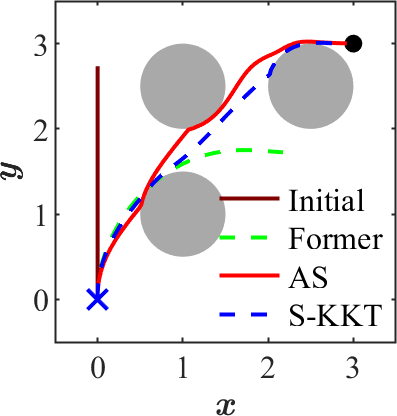}
  \subcaption{Initial trajectory is Feasible.}
 \end{minipage}
 \begin{minipage}[b]{0.31\linewidth}
  \centering
   \includegraphics[width=\columnwidth]{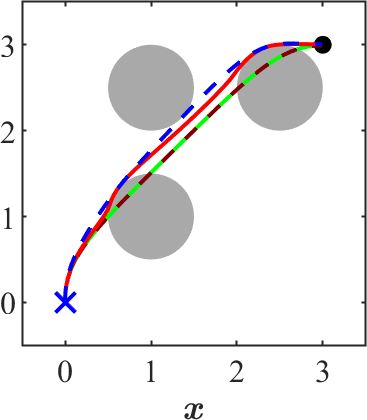} \subcaption{Slightly infeasible.}
 \end{minipage}
 \begin{minipage}[b]{0.31\linewidth}
  \centering
   \includegraphics[width=\columnwidth]{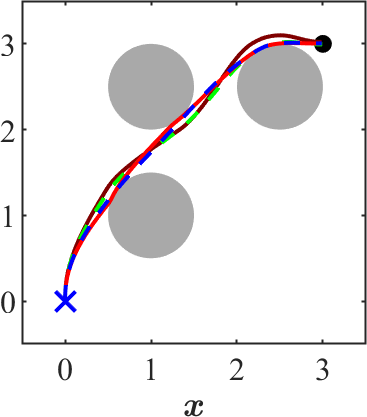} \subcaption{Close to optimal solution.}
 \end{minipage}
 \caption{2D car trajectories starting from several initial trajectories.}\label{fig:2dcar_init_trj}
\end{figure}

Fig. \ref{fig:2dcar_init_pnt} shows the result of the task starting from several initial points per algorithm. Optimization starts with six different initial points with no control.  
Fig. \ref{fig:2dcar_init_trj} shows the result of starting from initial trajectories. From left to right, a feasible trajectory, slightly infeasible trajectory, and close to the optimal trajectory were used as initial trajectories. \\
In the S-KKT algorithm, as we explained in eq. \ref{eq:s_init} and Algorithm \ref{alg:backpass}, $\bar{\bm{s}}$ is initialized by a positive value $\epsilon$.
This means that the algorithm regards the initial trajectory feasible even though it is not.
We experimentally confirmed that this is fine as long as the violation is small. In fact, this means S-KKT is able to handle initial trajectories that are slightly infeasible, which is something previous constrainted DDP algorithms cannot handle.

The procedure of the experiment is as follows.
We made several initial trajectories with different amounts of violations by changing the radii of obstacles and used them as initial trajectories of the original problem.
In our 2D car setting, original radii of obstacles are 0.5 m. We changed radius to 0.4 and ran the algorithm to get an optimized trajectory. Then, we used this trajectory as an initial trajectory with 0.1 violation for the original problem whose obstacles have 0.5 radii.
Our algorithm could successfully handle an initial violation up to 0.3.

\begin{figure*}[t!] 
\begin{subfigure}[b]{0.33\textwidth}
  \centering
  \includegraphics[width=\linewidth]{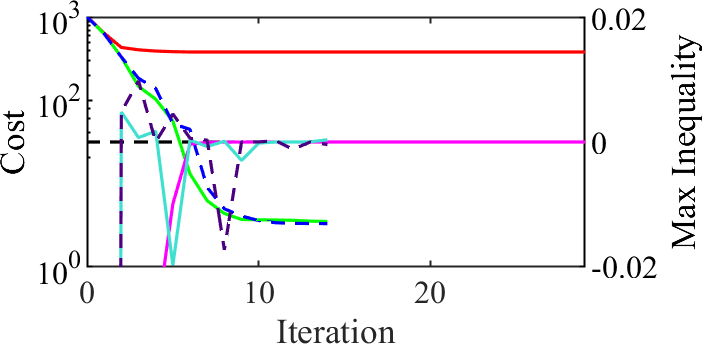} 
\end{subfigure}
\begin{subfigure}[b]{0.33\textwidth}
  \centering
  \includegraphics[width=\linewidth]{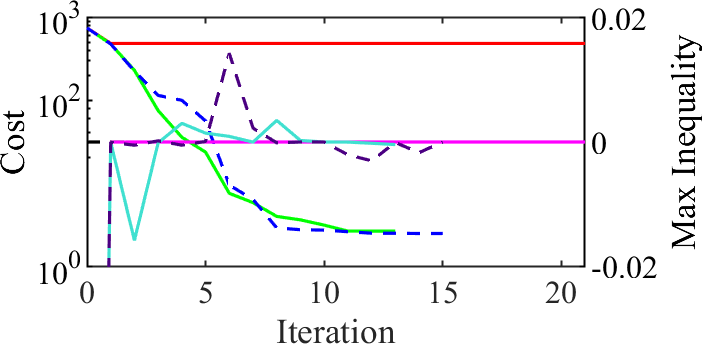} 
\end{subfigure}
\begin{subfigure}[b]{0.33\textwidth}
  \centering
  \includegraphics[width=\linewidth]{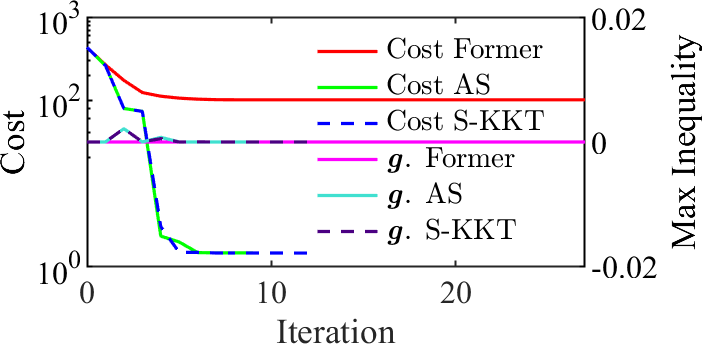} 
\end{subfigure}
\newline
\begin{subfigure}[b]{0.33\textwidth}
  \centering
  \includegraphics[width=\linewidth]{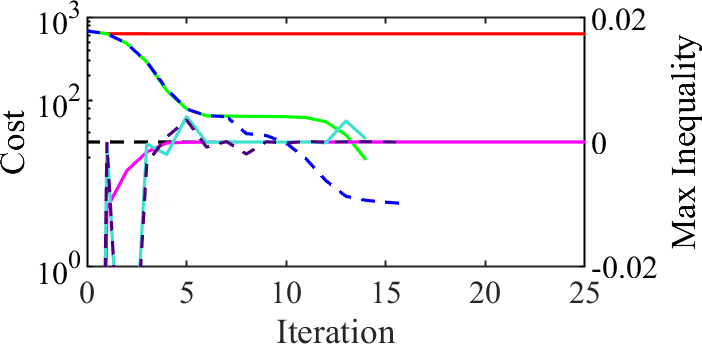} 
\end{subfigure}
\begin{subfigure}[b]{0.33\textwidth}
  \centering
  \includegraphics[width=\linewidth]{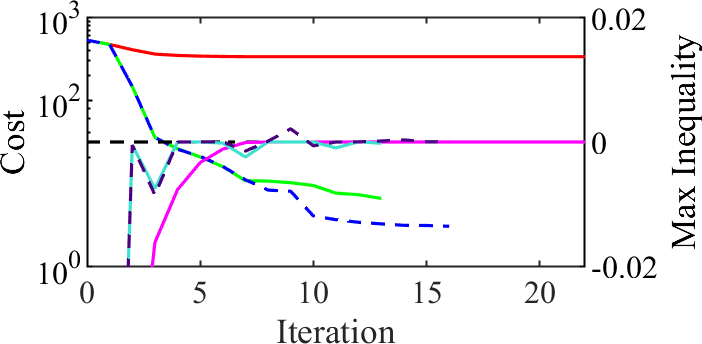}  
\end{subfigure}
\begin{subfigure}[b]{0.33\textwidth}
  \centering
  \includegraphics[width=\linewidth]{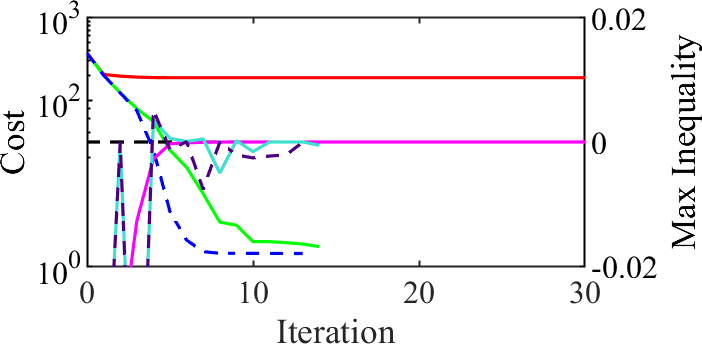}  
\end{subfigure}
\caption{Cost and max. inequality constraint of 2D car starting from several initial points.}
\label{fig:2dcar_init_pnt_cst}
\end{figure*}

\begin{figure*}[!]
\begin{subfigure}[b]{0.33\textwidth}
  \centering
  \includegraphics[width=\linewidth]{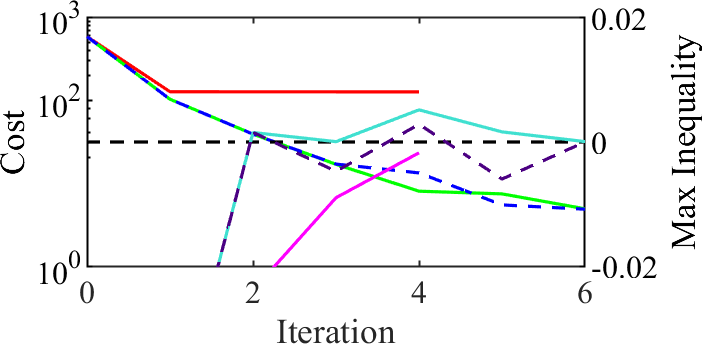}  
  \caption{Starting from feasible trajectory.}
\end{subfigure}
\begin{subfigure}[b]{0.33\textwidth}
  \centering
  \includegraphics[width=\linewidth]{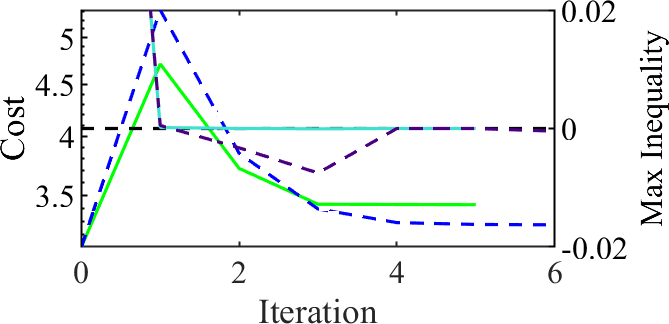}  
  \caption{Slightly infeasible trajectory.}
\end{subfigure}
\begin{subfigure}[b]{0.33\textwidth}
  \centering
  \includegraphics[width=\linewidth]{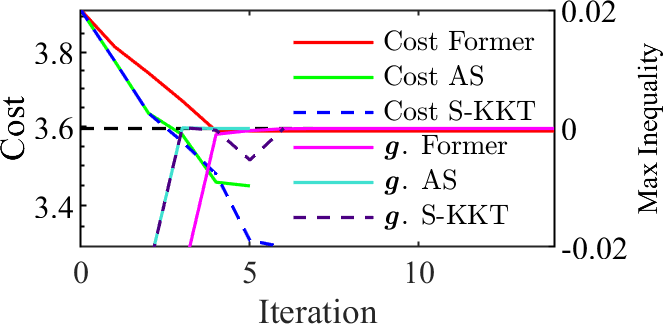}  
 \caption{Close to the optimal trajectory.}
\end{subfigure}
  \caption{Cost and max. inequality constraint of 2D car starting from several initial trajectories.}
  \label{fig:2dcar_init_trj_cst}
\end{figure*}

Fig. \ref{fig:2dcar_init_pnt_cst} shows the relationship between cost, max. value of constrained function $\bm{g}$ and iteration initialized with several start points. Fig. \ref{fig:2dcar_init_trj_cst} shows that of starting with initial trajectories. We specified the number of maximum outer iterations to be 15. The algorithm stops either when the max iteration is reached, when the regularizers are larger than the criteria, or the gradient of objective is very small. Introducing the slack variable in our method makes the trajectory smoother and we obtain the lowest converged cost. Our method could also get out of the prohibited region.


\subsubsection{Quadroter}
\begin{figure*}[!]
\begin{subfigure}[b]{0.33\textwidth}
  \centering
  \includegraphics[width=\linewidth]{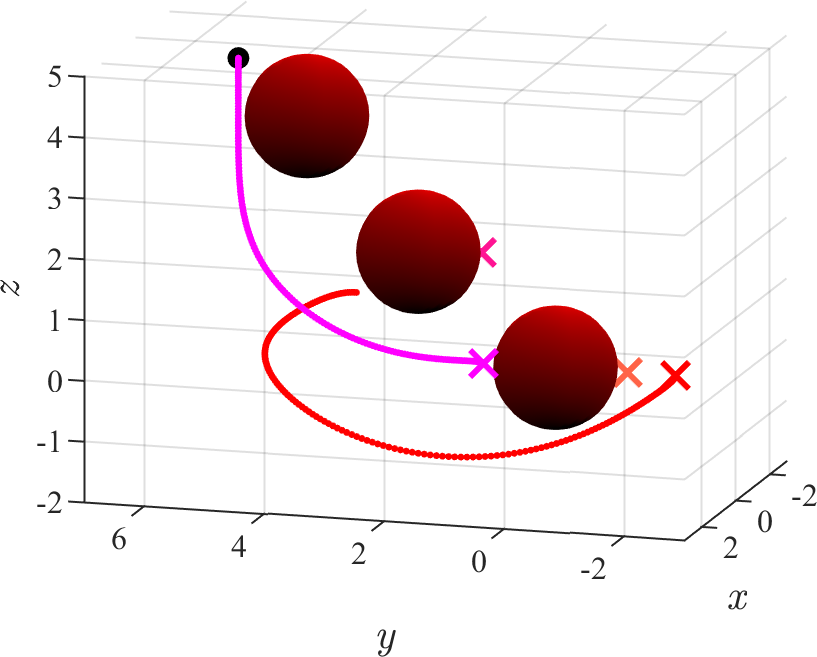}  
  \caption{Former algorithm}
  \label{fig:sub-first}
\end{subfigure}
\begin{subfigure}[b]{0.33\textwidth}
  \centering
  \includegraphics[width=\linewidth]{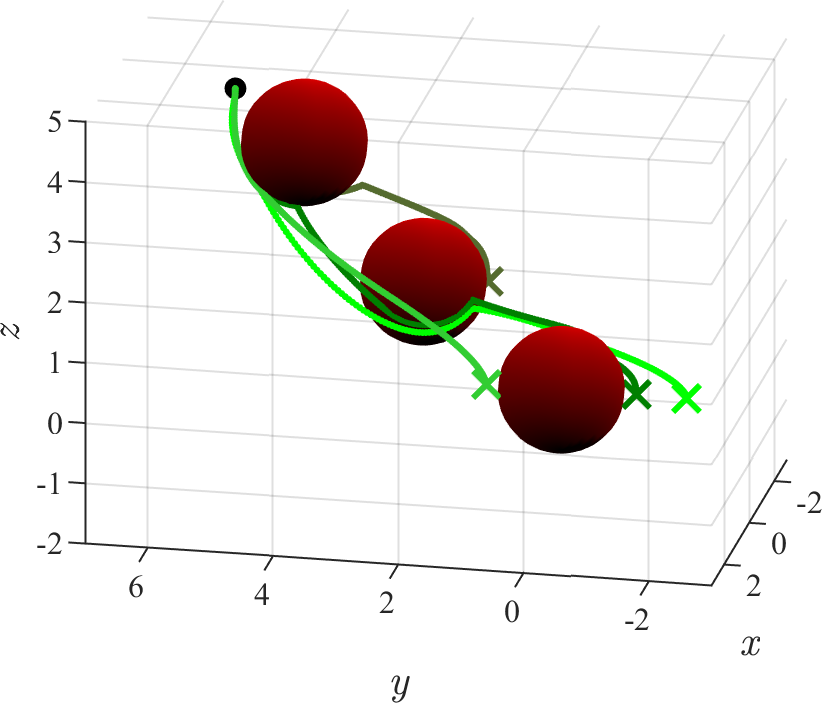}  
  \caption{Active set method}
  \label{fig:sub-second}
\end{subfigure}
\begin{subfigure}[b]{0.33\textwidth}
  \centering
  \includegraphics[width=\linewidth]{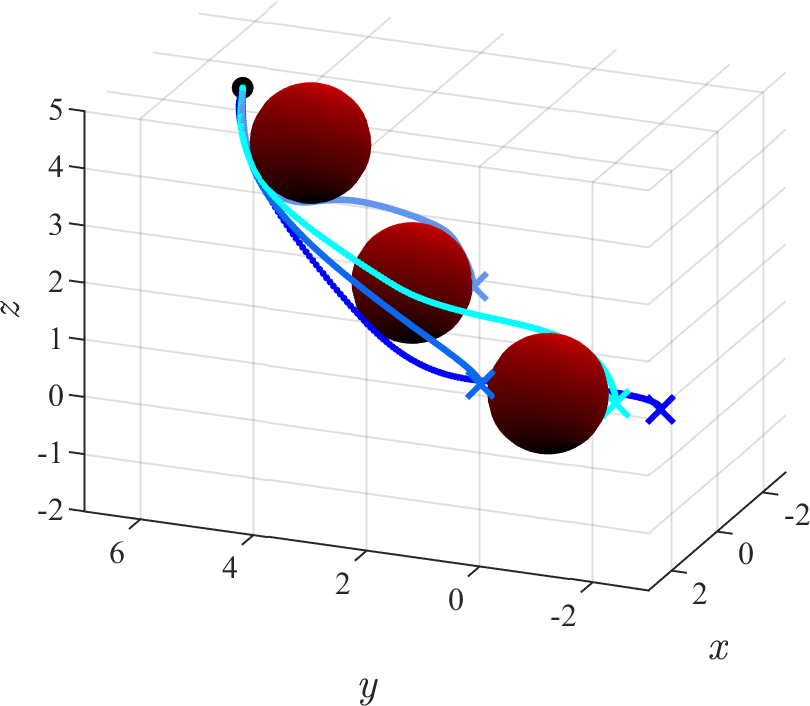}  
  \caption{S-KKT}
\end{subfigure}
\caption{Trajectories from different algorithms. Starting points are shown by ``x".}
  \label{fig:quad}
\end{figure*}

\begin{figure*}[!]
\begin{subfigure}[b]{0.5\textwidth}
  \centering
  \includegraphics[width=\linewidth]{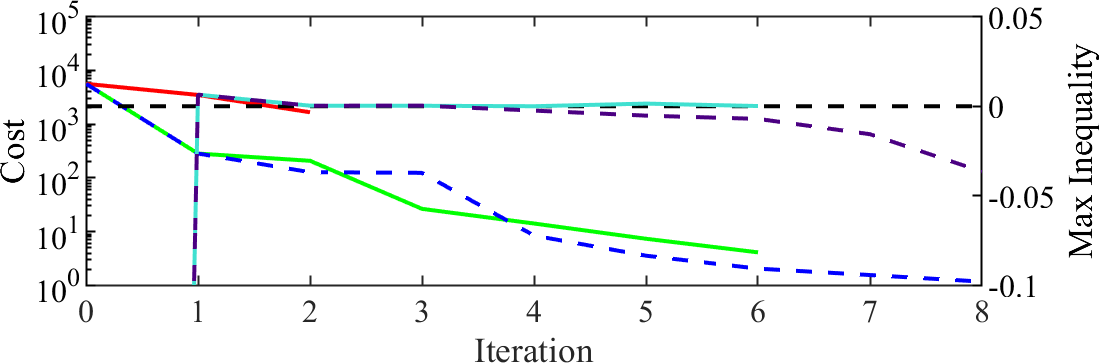}  
\end{subfigure}
\begin{subfigure}[b]{0.5\textwidth}
  \centering
  \includegraphics[width=\linewidth]{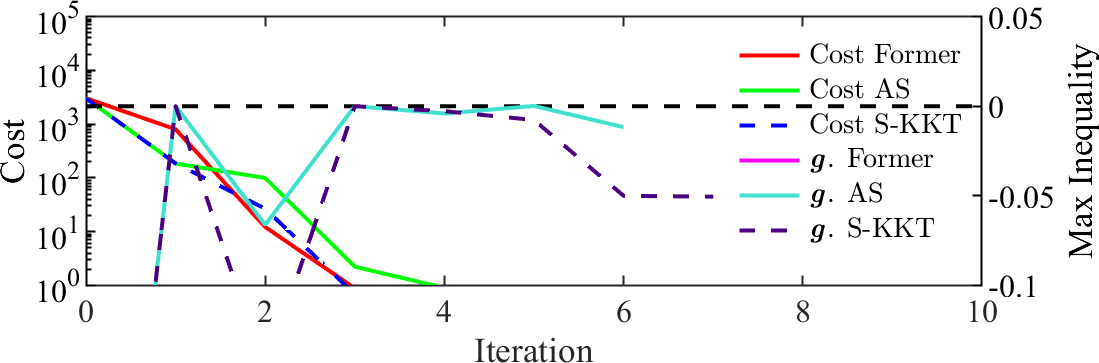}  
\end{subfigure}
\begin{subfigure}[b]{0.5\textwidth}
  \centering
  \includegraphics[width=\linewidth]{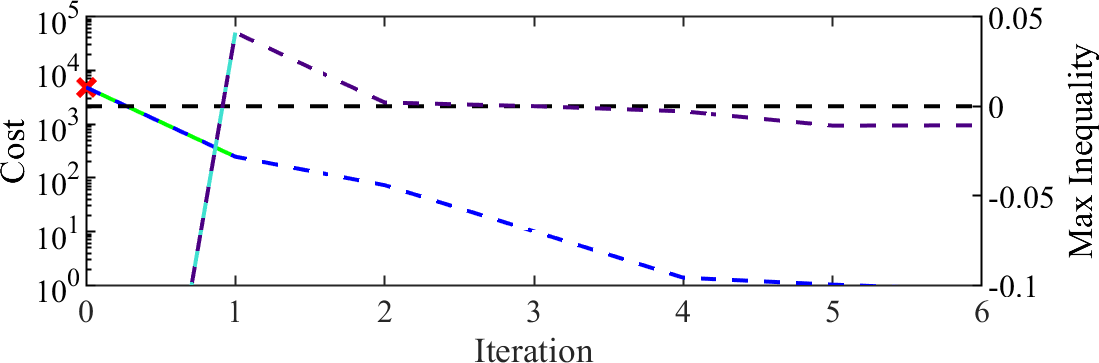}  
\end{subfigure}
\begin{subfigure}[b]{0.5\textwidth}
  \centering
  \includegraphics[width=\linewidth]{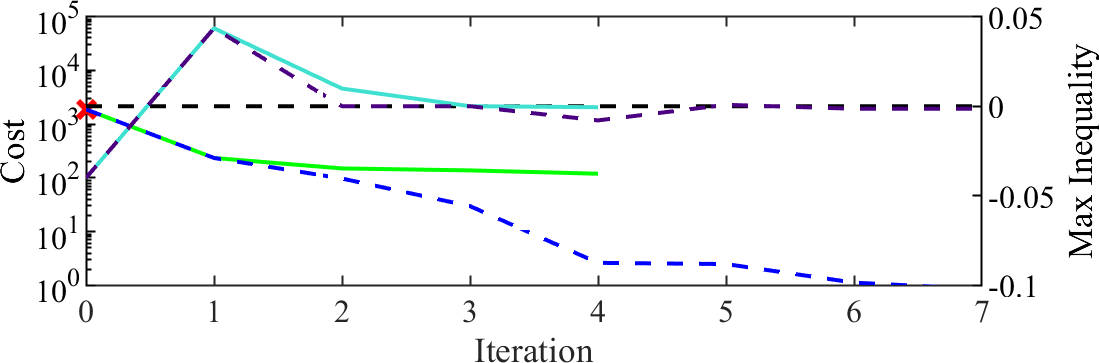}  
\end{subfigure}
\caption{Cost and max. inequality constraints of quadroter task.}
 \label{fig:quad_cost}
\end{figure*}

{\label{subsubsec:quad}}
We test our algorithm on a quadroter system \cite{Luukkonen2011}. The quadroter reaches a goal $\bm{x}_{g} = [1, 5, 5]^{\mathsf{T}}$ avoiding three obstacles. Fig. \ref{fig:quad} shows trajectories starting with four different initial hovering points from three different algorithms, that is, former KKT algorithm, active set method, and S-KKT. And Fig. \ref{fig:quad_cost} shows the cost and max. value of the constrained function $\bm{g}$. Again, S-KKT has the best performance.
    
\subsection{Combination of S-KKT and Augmented Lagrangian method}
\subsubsection{Control constraints}
\begin{figure}[h]
 \begin{minipage}[b]{0.49\linewidth}
  \centering
  \includegraphics[width=\columnwidth]{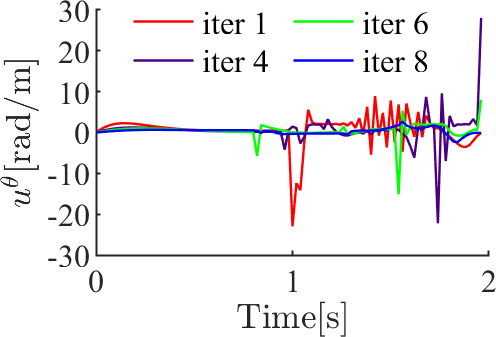}
  \subcaption{control $u^{\theta}$}
 \end{minipage}
 \begin{minipage}[b]{0.49\linewidth}
  \centering
  \includegraphics[width=\columnwidth]{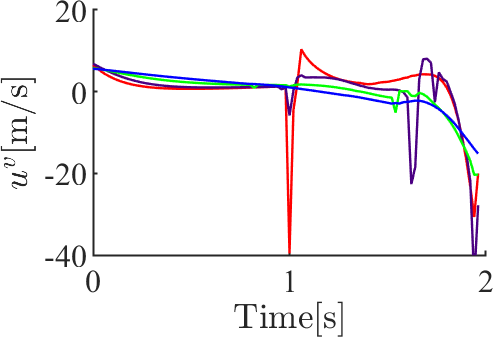}
 \subcaption{control $u^v$}
 \end{minipage}
 \caption{control over iteration}\label{fig:ctrl_iter}
\end{figure}

\begin{figure}
    \centering
    \includegraphics[scale = 0.5]{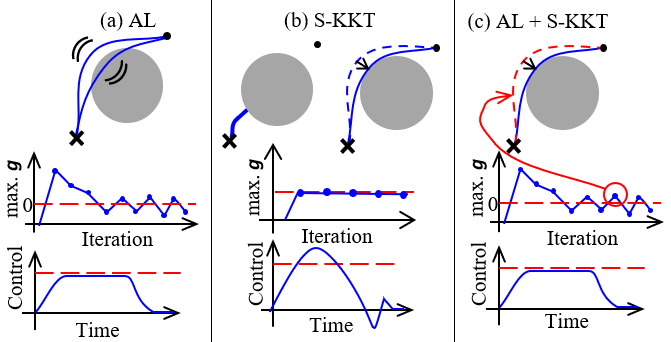}
    \caption{Concept of combining AL and S-KKT.\
    (a) AL-DDP shows oscillatory behavior around the boundary and take many iterations
  to keep feasibility.\
(b) S-KKT get stuck or infeasible when
large control is required. This situation happens when goal state is far form initial trajectory or large control is required to dodge the obstacle.\
Be able to keep feasibility.
(c) S-KKT is good at pushing 
 a trajectory which is
 close to the optimal.\
Obtain initial trajectory from AL and optimize the trajectory
  more to the optimal by S-KKT.}
    \label{fig:combination_concept}
    \end{figure}

Because S-KKT can take constraints in consideration only two time steps forward, sometimes state and control constraints conflict with each other.
In the 2D car case, for example, a car is trying to reach the target and suddenly finds an obstacle. If the control is not limited, the car can quickly steer to dodge the obstacle. However, if the control is limited, it cannot perform a sudden turn and makes a collision, making the trajectory infeasible. 
Fig. \ref{fig:ctrl_iter} shows how the control changes over iterations when the control is not limited. In this example, a 2D car starts from a static point [0, 0, 0, 0] with $\bm{0}$ initial control and reaches $\bm{x}_g$ explained in Section \ref{subsubsec:2Dcar}. In early iterations, large control spikes are observed. These spikes get smaller in future iterations, because the high control is penalized in the cost function. We can expect the optimizer to make the control spikes smaller than the arbitrary control limits, but there is no guarantee. Therefore S-KKT cannot explicitly apply control constraints to a trajectory as it is.
We solve this problem by combining AL-DDP and S-KKT. Using the control limited DDP technique in the backward pass, AL-DDP can apply control constraints to a trajectory. AL-DDP is very fast for the first few iterations but gets slow when it comes close to the boundary of constraints and sometimes violates the constraints.\\
Usually the trajectory oscillates around the boundary. Whereas S-KKT takes relatively a longer time for one iteration, but can keep feasibility. Though S-KKT can not handle a problem in which state and control constraints conflict each other. Typically the conflict happens when the initial state is far from the goal state and it needs to be changed a lot. However, given a good initial trajectory, S-KKT is good at pushing it close to boundary as shown in Fig. \ref{fig:2dcar_init_trj}. We feed S-KKT with the output of AL-DDP and optimize it under the expectation that large control is not required.
The concept of the combination is shown in Fig \ref{fig:combination_concept}. After receiving an optimized trajectory from AL-DDP, S-KKT solves the QP problem in its forward pass shown in 
eq. \eqref{forawrd_QP} with additional box control constraints,
\begin{align}
\bm{u}_{l} - \bar{\bm{u}}_k \leq  \delta \bm{u}_k \leq  {\bm{u}}_{u}-\bar{\bm{u}}_k,
\end{align}
If the QP problem is infeasible at time step $t_k$, we make the control more conservative by multiplying  $0 < \eta < 1$ to the box constraints as we do in \eqref{eq:trust_region},
\begin{align}
\eta(\bm{u}_{l}-\bar{\bm{u}}_k) \leq  \delta \bm{u}_k \leq  \eta({\bm{u}}_{u}-\bar{\bm{u}}_k),
\end{align}
and resolve the QP problem again from $t_0$ until the problem can be solved over the entire time horizon. This strategy, making a trajectory closer to a former one until it becomes feasible, is not good when the initial trajectory is far from desired one, and/or large control is required to dodge the obstacles as shown in Fig. \ref{fig:ctrl_iter}. In our case, however, thanks to a good initial trajectory from AL-DDP, this strategy fits well with S-KKT.\\
\subsubsection{2D car}
\begin{figure*}[!]
\begin{minipage}{.32\textwidth}
\begin{subfigure}{\textwidth}
\includegraphics[width=\textwidth]{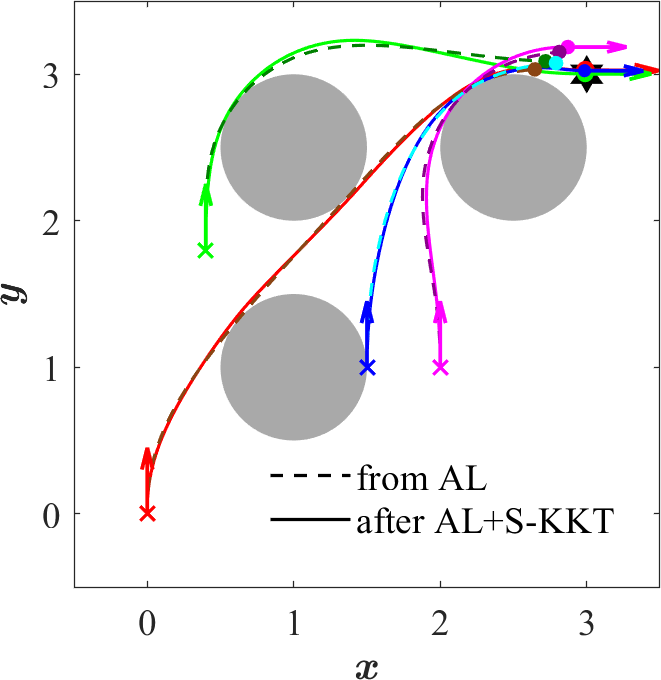}
\subcaption{Trajectories reaching a goal $[3,3,\frac{\pi}{2},0]^{\mathsf{T}}$ starting from static initial points.}
\label{fig:results_limited_trj}
\end{subfigure}
\end{minipage}
\hfill
\begin{minipage}{.333\textwidth}
\begin{subfigure}{\textwidth}
\includegraphics[width=\textwidth]{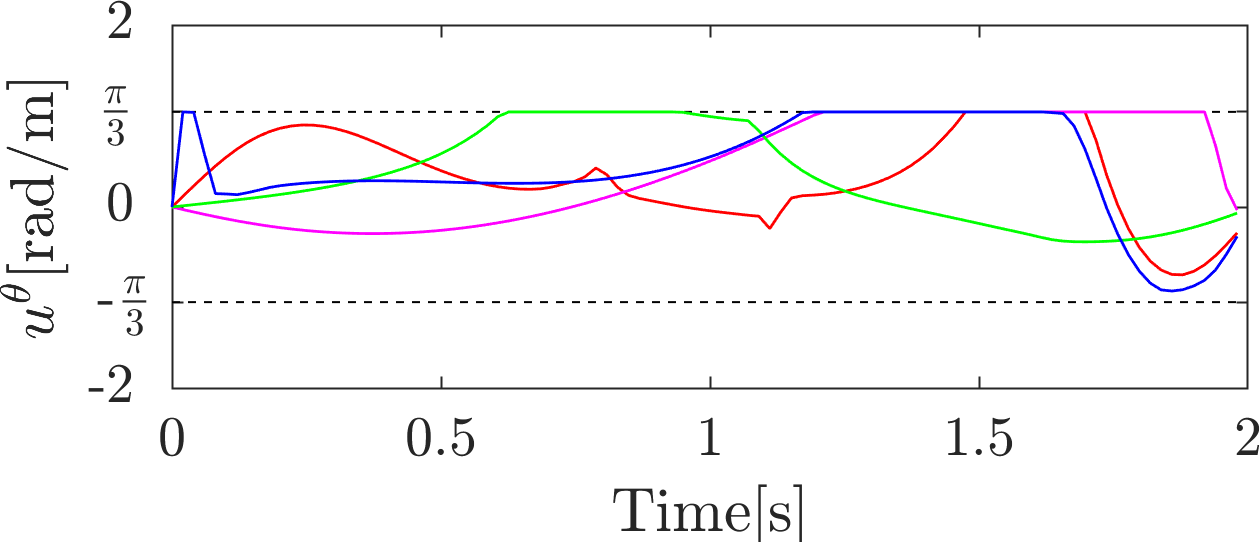}

\subcaption{Constrained $u^{\theta}$}
\label{fig:results_limited_ctrl_th}
\end{subfigure}
\begin{subfigure}{\textwidth}
\includegraphics[width=\textwidth]{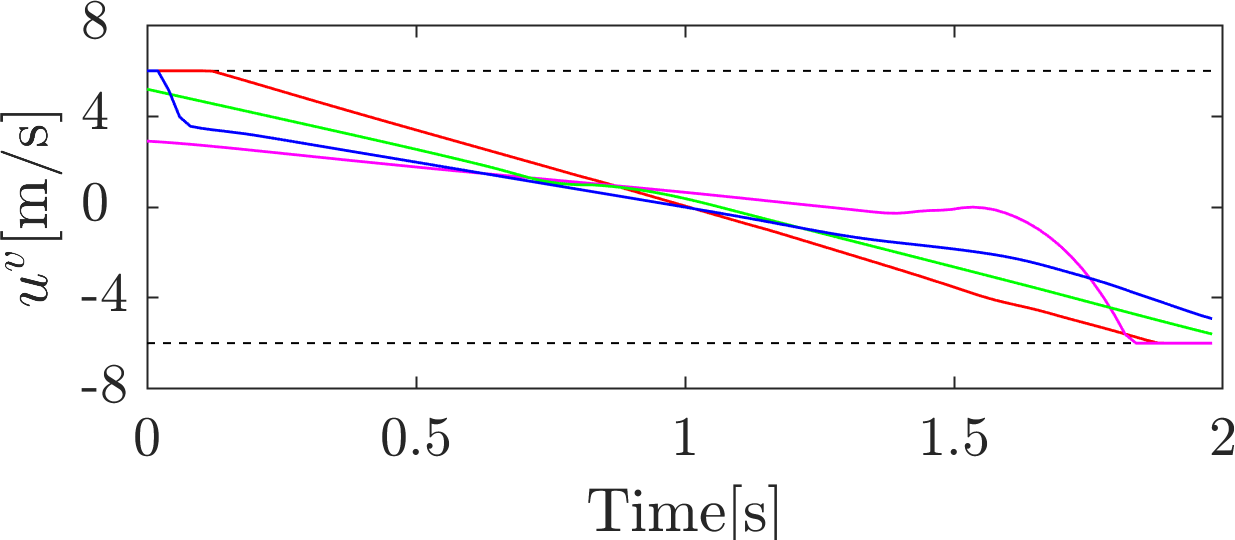}
\subcaption{Constrained $u^v$}
\end{subfigure}
\end{minipage}
\hfill
\begin{minipage}{.333\textwidth}
\begin{subfigure}{\textwidth}
\includegraphics[width=\textwidth]{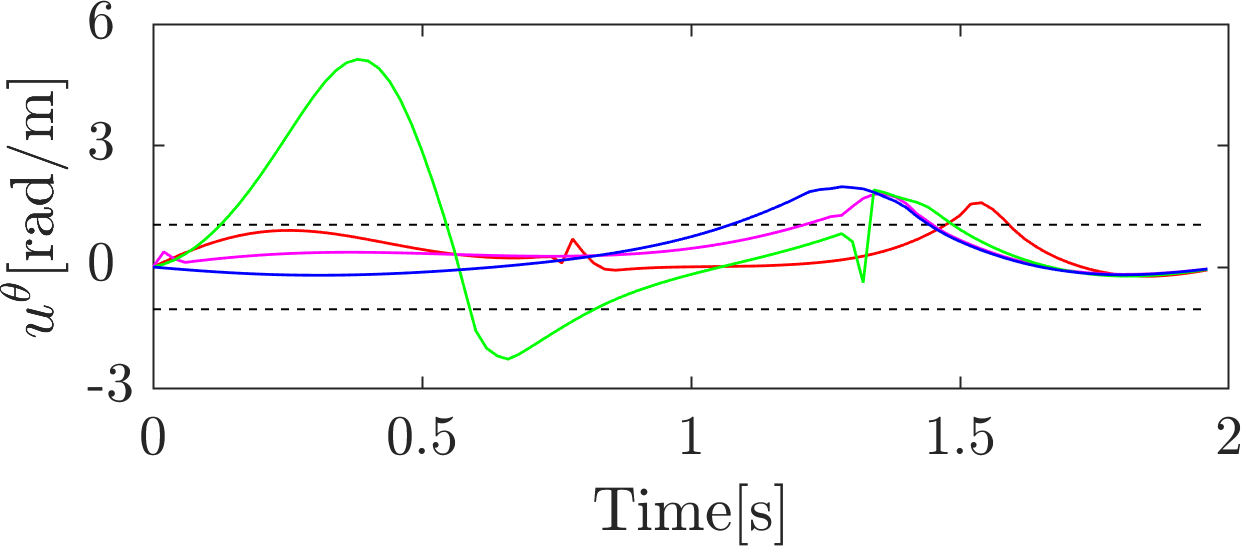}
\subcaption{Unconstrained $u^{\theta}$}
\label{fig:results_no_limit_ctrl_th}
\end{subfigure}
\begin{subfigure}{\textwidth}
\includegraphics[width=\textwidth]{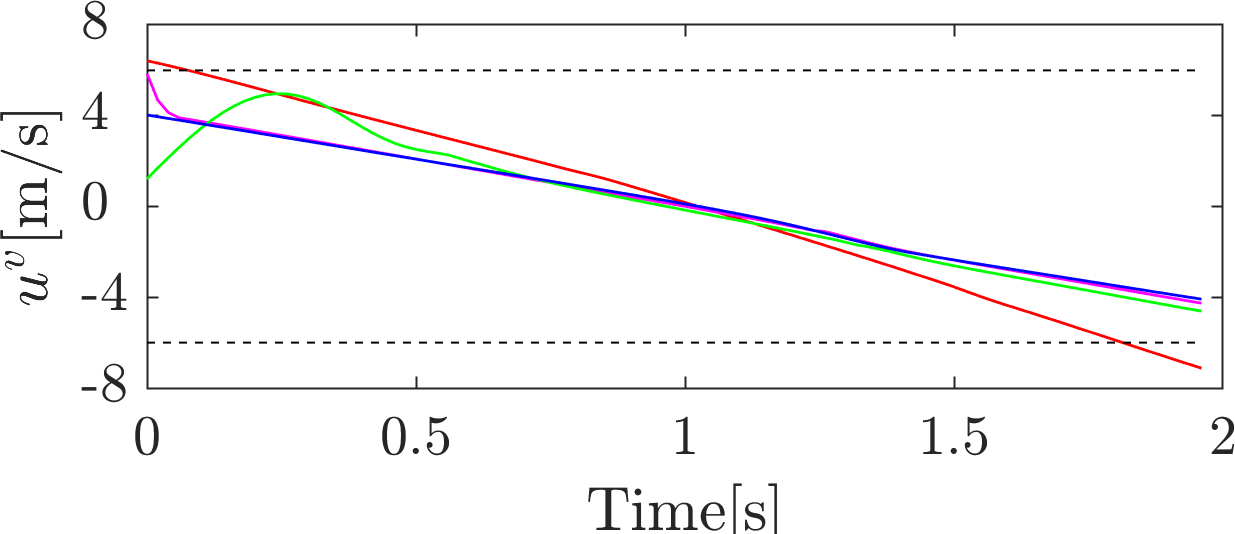}
\subcaption{Unconstrained $u^v$}
\label{fig:results_no_limit_ctrl_vel}
\end{subfigure}
\end{minipage}
\caption{Results from combined algorithm. In Fig. \ref{fig:results_limited_trj}, position and orientation at initial and final time steps from final result are shown by points and arrows. For AL, only final positions are plotted as circles, and arrows are omitted to make the figure easy to see. }  
\label{fig:results_limited}
\end{figure*}

To examine the performance of the combination,
we used the same problem setting of the 2D car in Section \ref{subsubsec:2Dcar}, and applied control limit as, 
\begin{align}
    -\frac{\pi}{3} \leq u^{\theta} \leq \frac{\pi}{3}, \quad 
    -6 \leq  u ^v \leq  6
\end{align}
The results and comparison between unconstrained control case are shown in Fig. \ref{fig:results_limited}. The algorithm could successfully handle both state and control constraints.

We observed that when the steering control constraint was too tight, the car could
only satisfy the desired angle (see pink trajectory in Fig.  \ref{fig:results_limited_trj}). In its control graph in Fig. \ref{fig:results_limited_ctrl_th}, we can see that maximum steering control was applied to dodge the obstacle and to reach the goal but it was not enough.
As shown in the green trajectory in Fig. \ref{fig:results_limited_trj}, it reached the goal, dodging the first obstacle from the top. Whereas in the control unconstrained case in Fig. \ref{fig:2dcar_init_pnt}, it could make a sharp turn and dodge the first obstacle from the bottom. This change can be also seen by comparing Fig. \ref{fig:results_limited_ctrl_th} and Fig. \ref{fig:results_no_limit_ctrl_th}. In the constrained case,

\subsection{Performance Analysis}
Next we perform a thorough analysis between five algorithms, SQP, S-KKT based DDP, AL-DDP, AL-DDP with SQP, and AL-DDP with S-KKT. The SQP algorithm used in this comparison is the one available through Matlab's Optimization Toolbox. For the rest of the algorithms, we have implemented them ourselves. In the five algorithms, the last two,  AL-DDP with SQP and AL-DDP with S-KKT are a combination of two different optimization algorithms. In these two combination methods, the algorithms start optimizing using AL-DDP and switch to SQP or S-KKT as explained in \ref{sec:ALandKKT}. We compare the five algorithms in terms of performance metrics, namely cost, time, and feasibility, in three systems, cart pole, 2D car, and quadroter. We also specify different time horizons for the same examples. Here the time horizon $H$ is the number of time steps and the time step size $dt$ is fixed in each example.
Furthermore, we perform experiments with a time budget as well as letting the algorithms run until they have reached the convergence criteria. All of the simulations are performed in Matlab 2019b on a CPU architecture of i7-4820K (3.7 GHz).

Note we divide feasibility of the solution into two parts, one is feasibility with respect to the constraint function $\bm{g}$ and the other is with respect to the dynamics $\bm{f}$. SQP handles them as inequalities (for $\bm{g}$) and equalities (for $\bm{f}$) constraints, whereas DDP-based methods can implicitly satisfy dynamics since they are used during the optimization process.

\subsubsection{Exit criteria}
The exit condition of optimization was based on two criteria as shown in Table \ref{tab:exit_criteria_combo}. One was constraint satisfaction, where $ 1 \times 10^{-7} $ was used for all the algorithms. The other was an optimality criterion. For the DDP based methods, we used the change in the optimization objective between iterations, set to $8 \times 10^{-2}$. For SQP, we used an optimality condition shown in the first equation of \eqref{eq:first_optimality}. This condition was set by choosing the Matlab \texttt{fmincon} option {\texttt{OptimalityTolerance}} to be $ 1 \times 10^{-2} $. We first set the value to be $8 \times 10^{-2}$ which is the same as that of other DDP based methods. However, SQP stopped its optimization process reaching local minima in the early iterations for several examples, which kept the cost of SQP relatively much higher. In the 2D car and the quadroter case for example, the cost was about a hundred times higher than S-KKT. Therefore, we made the criteria smaller to further the SQP optimization process. When using the change of objective function instead of the optimality tolerance in SQP, we observed that the condition was also satisfied in the early iterations of the optimization, but with a large constraint violation. This means that only the constraint satisfaction was working effectively as the exit condition. Thus, we decided to use the optimality tolerance for SQP.
For AL-SQP and AL-S-KKT, we also have conditions on exiting the AL optimization scheme and switching to the next scheme as shown in Table \ref{tab:exit_criteria_combo_b}). In lieu of fairness, we decided to keep this "switching condition" and other AL parameters the same between both algorithms even if they may have benefited with different conditions for the overall convergence requirements. For these AL schemes, the constraint satisfaction tolerance was $ 1 \times 10^{-2}$ and the bound on the change of the cost was $1$.

\begin{table}
\caption{Exit criteria}\label{tab:exit_criteria_combo}
\begin{subtable}{1\columnwidth}
\subcaption{Exit criteria for single algorithms.}
\label{tab:exit_criteria_single}
\centering
   \begin{tabular}{c|ccc}
       criterion & SQP & S-KKT & AL-DDP\\ 
       \hline
      cnst. satisfaction & 1E-07 & 1E-07 & 1E-07\\
      \hline
      closeness to & opt. tolerance & change in cost & change in cost\\
      optimal solution& 1E-2 & 8E-2 &8E-2
   \end{tabular}
\end{subtable}
\newline
\vspace*{5 mm}
\newline
\begin{subtable}{1\columnwidth}
      \subcaption{Exit criteria for combination algorithms.}
   \label{tab:exit_criteria_combo_b}
\centering
   {\begin{tabular}{c|cc}
      stage & AL-SKKT & AL-DDP\\ 
      \hline
      first & \multicolumn{2}{c}{cnst. satisfaction: 1E-2}\\
      AL-DDP & \multicolumn{2}{c}{change in cost: 1}\\
      \hline
      final & \multicolumn{2}{c}{cnst. satisfaction:1E-7}\\
      optimization &opt. tolerance: 1E-2  &change in cost: 8E-2
   \end{tabular}}
\end{subtable}
\end{table}

\subsubsection{Cart Pole}
Table \ref{tab:comp_cpole} shows the results of the simulations for balancing a cart pole system. 
The system has four dimension of state $\bm{x}$, that is position of the cart $x$, its velocity $\dot{x}$, angle of the pendulum $\theta$, and angular velocity $\dot{\theta}$. The control of the system is thrust force $u$ applied to the cart. The dynamics is given as follows:  
\begin{align}
    \ddot{x}_{k} &= \frac{u_{k}-b\dot{x}_{k}+m(l+{\dot{\theta}_{k}}^2 -g\cos{\theta_{k}})\sin{\theta_{k}}}{(M+m\sin^2{\theta_{k}})},\\\notag
    \ddot{\theta}_{k} &= \frac{g(M+m)\sin{\theta_{k}}-(u-b\dot{x}_{k}+ml\dot{\theta}_{k}^{2}\sin{\theta_{k}})\cos{\theta_{k}}}{l(M+m\sin^2{\theta_{k}})},
\end{align}
where $M$ is a mass of the cart, $m$ is that of pendulum, $l$ is a length of the arm of the pendulum, $g$ is gravitational acceleration,  and $b$ is a coefficient of friction between the car and the floor. The problem has constraints in the position and angle:
\begin{align}
\bm{g}_{\rm{cp}} =
\begin{bmatrix}
x^2 - x_{\rm{lim}}^2\\
\theta - \theta_{\rm{lim}}
\end{bmatrix} \leq \bm{0}
\end{align}
\begin{table}
    {\fontsize{6.7}{10}\selectfont
    \centering\captionsetup{justification = centering}
    \caption{Performance metrics for a cart pole system until convergence.}
    \label{tab:comp_cpole}
    \begin{tabular}{c|c|ccccc}
        $H$ & Metric & SQP & S-KKT & AL & AL-SQP & AL-S-KKT \\
        \hline
        & Cost & 3.11 & 5.90 & 5.97 & $\bm{3.06}$ & 6.06 \\
        & Time & 12.3 & 7.38 & 5.32 & 11.7 & $\bm{4.70}$ \\
        100 & Feas. ($g$) & 0 & 0 & 1.89E-04 & 0 & 0 \\
        & Feas. ($f$) & 1.90E-03 & 0 & 0 & 3.19E-04 & 0 \\
        
        \hline
        & Cost & 2.80 & 7.77 & 5.81 & $\bm{2.77}$ & 5.28 \\
        & Time & 52.2 & $\bm{12.0}$ & 10.7 & 48.5 & 12.1 \\
        200 & Feas. ($g$) & 0 & 0 & 2.50E-04 & 0 & 0 \\
        & Feas. ($f$) & 2.78E-05 & 0 & 0 & 1.98E-04 & 0 \\
        \hline  
    \end{tabular}
    }
\end{table}
\begin{table}[!]
    {\fontsize{6.7}{10}\selectfont
    \centering\captionsetup{justification = centering}
    \caption{Performance metrics for a cart pole system with time budget.}
    \label{tab:comp_cpole_budget}
    \begin{tabular}{c|c|ccccc}
        $H$ & Metric & SQP & S-KKT & AL & AL-SQP & AL-S-KKT \\
        \hline
        & Cost & 3.53 & 5.90 & 5.97 & $\bm{3.18}$ & 6.06 \\
        & Time & 6 & 6 & 5.25 & 6 & $\bm{4.69}$ \\
        100 & Feas. ($g$) & 0 & 0 & 1.89E-04 & 0 & 0 \\
        & Feas. ($f$) & 4.12E-03 & 0 & 0 & 4.83E-03 & 0 \\
        \hline
        & Cost & 7.21 & 8.37 & 7.33 & 5.80 & $\bm{7.14}$ \\
        & Time & 6 & 6 & 6 & 6 & 6 \\
        200 & Feas. ($g$) & 0 & 0 & 3.34E-04 & 5.50E-06 & 0 \\
        & Feas. ($f$) & 1.22E-03 & 3.55E-06 & 0 & 3.51E-02 & 0 \\
        \hline  
    \end{tabular}
    }
\end{table}
Pure SQP performed the slowest with dynamics violation, although it achieved a very low cost. It required a much longer time for a longer time horizon compared to other methods. This is understandable because in SQP, a longer time horizon corresponds to a larger matrix (which needs to be inverted) containing the equality constraints for the dynamics.
S-KKT also takes time and it accrues a high cost compared to SQP. However, it does satisfy feasibility. AL-DDP on its own cannot reach the same levels of constraint satisfaction as S-KKT, which makes sense since the AL approach oscillates near the constraint bounds, but converges faster. When pairing AL with SQP, there is no significant change compared to original SQP. Pairing AL with S-KKT, however, we see an improvement. In the case of $H = 100$, compared to S-KKT, AL-S-KKT converges faster to an almost equally low cost. In the case of $H = 200$, AL-S-KKT takes slightly longer time, but converges to a lower cost.
Compared to AL, AL-S-KKT takes a long time, but satisfies feasibility.
From Table \ref{tab:comp_cpole_budget}, we can see the longer time horizon decreases the performance of SQP in terms of speed and constraint satisfaction, where AL-S-KKT is not affected as much. In S-KKT, there is a small violation of constraint possibly from the linearization error of the constraint function. The solution may satisfy the linearized constraint in \eqref{eq:g_expanded}, but not the original one. The error decreases as $\delta{\bm{u}}$ decreases, but if we use a time budget and stop the optimization process before convergence, there is a possibility that the solution has a small violation.

\subsubsection{2D car}
We use the same problem setting as \ref{subsubsec:2Dcar}. For these metrics we initialize the problem with six different starting points and take the average of the cost, time, and feasibility. The time step used was $dt = 0.02$ s. Table \ref{tab:comp_car} shows the result when we let run the algorithm until convergence, and Table \ref{tab:comp_car_budget} shows the result under time budget. In the case of $H = 100$, time budget was 3 s, and when $H = 200$, it was 6 s. In this example algorithms behave similarly as the example of a cart pole, and combination methods show their performance more clearly in terms of speed.  
\begin{table}[h]
    {\fontsize{6.7}{10}\selectfont
    \centering\captionsetup{justification = centering}
    \caption{Performance metrics for a 2D car system until convergence.}
    \label{tab:comp_car}
    \begin{tabular}{c|c|ccccc}
        $H$& Metric & SQP & S-KKT & AL & AL-SQP & AL-S-KKT \\
        \hline
        & Cost & 2.85 & 2.58 & 2.55 & $\bm{2.39}$ & 2.45 \\
        & Time & 11.1 & 5.51 & 5.69 & 5.77 & $\bm{2.41}$ \\
        100 & Feas. $g$ & 0 & 0 & 4.05E-04 & 0 & 0 \\
        & Feas. $f$ & 1.30E-06 & 0 & 0 & 5.01E-08 & 0 \\
        
        \hline
        & Cost & 1.66 & 1.20 & 1.12 & $\bm{0.995}$ & 1.05 \\
        & Time & 67.3 & 11.6 & 10.6 & 38.0 & $\bm{3.98}$ \\
        200 & Feas. $g$ & 0 & 0 & 1.26E-04 & -7.63E-08 & 0 \\
        & Feas. $f$ & 2.26E-05 & 0 & 0 & 4.63E-08 & 0 \\
        \hline  
    \end{tabular}
    }
\end{table}
\begin{table}[h]
    {\fontsize{6.7}{10}\selectfont
    \centering\captionsetup{justification = centering}
    \caption{Performance metrics for a 2D car system with time budget.}
    \label{tab:comp_car_budget}
    \begin{tabular}{c|c|ccccc}
        $H$ & Metric & SQP & S-KKT & AL & AL-SQP & AL-S-KKT \\
        \hline
        & Cost & 25.3 & 22.5 & 3.13 & 2.49 & $\bm{2.44}$ \\
        & Time & 3 & 2.77 & 3 & 3 & $\bm{2.36}$ \\
        100 & Feas. $g$ & 0 & 0 & 0 & 0 & 0 \\
        & Feas. $f$ & 1.02E-03 & 0 & 0 & 2.53E-05 & 0 \\
        \hline
        & Cost & 135 & 43.2 & $\bm{1.71}$ & 1.80  & 1.92 \\
        & Time & 6 & 6 & 6 & 6 & $\bm{3.88}$ \\
        200 & Feas. $g$ & 0 & 0 & 1.73E-03 & 0 & 0 \\
        & Feas. $f$ & 1.22E-02 & 0 & 0 & 1.10E-05 & 0 \\
        \hline  
    \end{tabular}
    }
\end{table}
\subsubsection{Quadroter}
In this example, we used same problem setting as Section \ref{subsubsec:quad} initialized with four different static hovering points, and take the average of performance metrics as we did in the 2D car example. We set a time step of $dt = 0.01$. As we can see from the results shown in Table \ref{tab:comp_quad}, SQP suffers from increase of dimension of the problem resulting in a much longer computational time. AL in the case of $H = 300$, could not get out from its inner optimization loop, and could not converge. We filled the corresponding table with ``N/A". Our S-KKT and AL-S-KKT, however could keep its stability and feasibility. In addition, they achieved a low cost in a short time. Table \ref{tab:comp_quad_budget} shows the result from the same problem under a time budget. For $H = 200$, the time budget was 6 s and for $H = 300$ it was 10 s. Single SQP took such a long time that it could not perform one single iteration, resulting in very high cost.
We have observed that our AL-S-KKT lost its speed performance affected by the first AL optimization process. AL consumed most of the time budget, allowing S-KKT only one or two iterations. We believe that more investigation or tuning of AL will make our AL-S-KKT much better.
\begin{table}[h]
    {\fontsize{6.67}{10}\selectfont
    \centering\captionsetup{justification = centering}
    \caption{Performance metrics for a quadroter system until convergence.}
    \label{tab:comp_quad}
    \begin{tabular}{c|c|ccccc}
        $H$& Metric & SQP & S-KKT & AL & AL-SQP & AL-S-KKT \\
        \hline
        & Cost & 5.87 & 7.95 & 8.10 & $\bm{5.5}$ & 7.03 \\
        & Time & 748 & 12.3 & 23.8 & 790 & $\bm{10.7}$ \\
        200 & Feas. $g$ & 0 & 0 & 0 & 0 & 0 \\
        & Feas. $f$ & 1.16E-05 & 0 & 0 & 2.84E-05 & 0 \\
        \hline
        & Cost & 5.91 & 5.99 & N/A & $\bm{5.37}$ & 6.88 \\
        & Time & 2.55E03 & 21.9 & N/A & 2.42E03 & $\bm{14.9}$ \\
        300 & Feas. $g$ & 0 & 0 & N/A & 0 & 0 \\
        & Feas. $f$ & 1.28E-05 & 0 & N/A & 2.08E-05 & 0 \\
        \hline  
    \end{tabular}
    }
\end{table}
\begin{table}[h]
    {\fontsize{6.8}{10}\selectfont
    \centering\captionsetup{justification = centering}
    \caption{Performance metrics for a quadroter system with time budget.}
    \label{tab:comp_quad_budget}
    \begin{tabular}{c|c|ccccc}
        $H$ & Metric & SQP & S-KKT & AL & AL-SQP & AL-S-KKT \\
        \hline
        & Cost & 2.43E03 & 26.1 & 10.7 & 8.12 & $\bm{7.48}$ \\
        & Time & 6 & $\bm{5.68}$ & 6 & 6 &5.91 \\
        200 & Feas. $g$ & 0 & 0 & 0 & 0 & 0 \\
        & Feas. $f$ & 1.86E-03 & 0 & 0 & 7.99E-04 & 0 \\
        \hline
        & Cost & 2.86E03 & $\bm{7.52}$ & 7.91 & 8.01  & 7.78 \\
        & Time & 10 & 10 & 10 & 10 & 10 \\
        300 & Feas. $g$ & 0 & 0 & 0 & 0 & 0 \\
        & Feas. $f$ & 2.50E-03 & 0 & 0 & 6.85E-09 & 0 \\
        \hline  
    \end{tabular}
    }
\end{table}
\section{Conclusion} 
\label{sec:conclusion}

 In this paper we have introduced novel constrained trajectory optimization methods that outperform previous versions of constrained  DDP. Some key ideas in this paper rely on the combination of slack variables together with augmented Lagrangian method and the KKT conditions. In particular,  
 
 \begin{itemize}
    \item Slack variables are an  effective way to get lower cost with respect to alternative algorithms relying on the  active set method. 
     \item The S-KKT method is able to handle both state and control constraints in cases where the feasibility set is small.
     \item The S-KKT methods is more robust to initial conditions of the state trajectory. 
     \item AL is very fast for first few iterations but get slow when it comes close to constraints and sometimes violate constraints. Whereas S-KKT takes time in one iteration, but can keep feasibility in a few iterations.
     By combining them we may be able to compensate for weakness of both and have a better algorithm. 
     
 \end{itemize}

 Future directions will include  mechanisms for uncertainty representations and learning, and development of  chance constrained trajectory optimization algorithms that have the benefits of the fast convergence of the proposed algorithms.



%
\appendices
%
\section{On the Augmented Lagrangian}
The penalty functions in \eqref{augmentedLagrangian} must be such that $\mathcal{P}'(y, \lambda, \mu):=\frac{\partial}{\partial y}\mathcal{P}(y,\lambda,\mu)$ is continuous for all $y\in\mathbb{R}$, $\lambda, \mu\in\mathbb{R}_{++}$ and: (i) $\mathcal{P}'(y, \lambda, \mu)\geq 0$, (ii) $\lim_{k\rightarrow\infty}\mu_{(k)}=\infty$ and $\lim_{k\rightarrow\infty}\lambda_{(k)}=\lambda>0$ imply that $\lim_{k\rightarrow\infty}\mathcal{P}'(y_{(k)}, \lambda_{(k)}, \mu_{(k)})=\infty$, (iii) $\lim_{k\rightarrow\infty}\mu_{(k)}=\infty$ and $\lim_{k\rightarrow\infty}\lambda_{(k)}=\lambda<0$ imply that $\lim_{k\rightarrow\infty}\mathcal{P}'(y_{(k)}, \lambda_{(k)}, \mu_{(k)})=0$.


\bibliographystyle{plainnat}
\bibliography{references}

\end{document}